\documentclass[a4paper,11pt]{amsart}
\usepackage{amssymb,amsmath,amsthm}
\usepackage[colorlinks=true,linkcolor=blue,citecolor=blue]{hyperref}
\hypersetup{pdftitle={A half-line of Schur idempotent norms},pdfauthor={Rupert H. Levene}}
\usepackage{graphicx}
\usepackage{tikz}
\usepackage{enumitem}
\DeclareMathOperator{\spn}{span}
\newtheorem{theorem}{Theorem}[section]
\newtheorem{lemma}[theorem]{Lemma}

\newtheorem{proposition}[theorem]{Proposition}
\theoremstyle{definition}
\newtheorem{definition}[theorem]{Definition}

\theoremstyle{remark}
\newtheorem{remark}[theorem]{Remark}
\newtheorem{question}[theorem]{Question}
\numberwithin{equation}{section}

\newcommand{\B}{\mathcal{B}}
\newcommand{\N}{\mathcal{N}}
\newcommand{\bN}{\mathbb{N}}
\newcommand{\bR}{\mathbb{R}}
\newcommand{\bZ}{\mathbb{Z}}
\newcommand{\sn}[1]{\|#1\|_\bullet}
\newcommand{\schur}{\bullet}

\newcommand{\Aut}{\operatorname{Aut}}

\title{A half-line of Schur idempotent norms}
\author{Rupert H. Levene}
\address{School of Mathematics and Statistics, University College Dublin, Belfield, Dublin 4, Ireland}
\email{rupert.levene@ucd.ie}

\begin{document}
\begin{abstract}
Using a family of bipartite graphs constructed by hanging identical rooted trees at every vertex of a highly symmetric graph, we show that every sufficiently large real number is the norm of an idempotent Schur multiplier on the bounded operators on $\ell^2$. This answers a question previously raised by the author.
\end{abstract}
\maketitle

\section{Introduction}

For index sets $I$ and $J$ and a complex matrix $B=[b_{ij}]_{i\in I,j\in J}$, the \emph{Schur norm} of $B$ is
\[\sn B=\sup\{\|B\schur X\|:X\in\B(\ell^2_J,\ell^2_I),\ \|X\|\le1\},\]
where $\|\cdot\|$ is the operator norm, $\ell^2_J$ is the Hilbert space $\ell^2$ of the index set $J$, the space $\B(\ell^2_J,\ell^2_I)$ consists of the bounded operators from $\ell^2_J$ to $\ell^2_I$, and $B\schur X=[b_{ij}x_{ij}]$ is the entrywise product of $B$ and $X=[x_{ij}]$. We call $B$ a \emph{Schur multiplier} if $\sn B<\infty$, and we say that $B$ is \emph{$0$--$1$} if $b_{ij}\in\{0,1\}$ for all $i$ and $j$.

For a $0$--$1$ Schur multiplier $B$, the map
\[\B(\ell^2_J,\ell^2_I)\to\B(\ell^2_J,\ell^2_I),\qquad X\mapsto B\schur X\]
is idempotent, and every idempotent Schur multiplier arises in this way. A \emph{bipartite graph} is a triple $G=(I,J,E)$ of disjoint sets $I$ and $J$, whose elements are the \emph{row vertices} and \emph{column vertices}, together with a set $E\subseteq I\times J$ of \emph{edges}. Its vertex set is $V(G)=I\sqcup J$. We identify a $0$--$1$ matrix $B$ indexed by $I\times J$ with the bipartite graph $G=(I,J,E)$ whose edges are the pairs $ij$ with $b_{ij}=1$, and we write $\sn G=\sn B$.

Let $\N$ be the set of norms of the idempotent Schur multipliers on $\B(\ell^2)$, that is,
\[\N=\{\sn B: B\text{ is a $0$--$1$ Schur multiplier indexed by }\bN\times\bN\}.\]
Livshits~\cite{Liv} showed that $\N\cap(1,\tfrac2{\sqrt3})=\emptyset$. Katavolos and Paulsen~\cite[Remark~5]{KatP} recorded an observation of Davidson that $\N$ is closed under multiplication and under suprema of bounded sets, and noted that $2$ is a limit point of $\N$ by~\cite{BCD}. The set $\N$ is also unbounded~\cite{KP}.

In~\cite[Theorem~1.1]{L14}, the author extended this result of Livshits. With $\eta_6=\frac25\sqrt{5+2\sqrt5}\approx1.231$, that theorem states that
\[\N\cap[0,\eta_6]=\Bigl\{0,\ 1,\ \sqrt{\tfrac43},\ \tfrac12(1+\sqrt2),\ \tfrac1{15}\sqrt{169+38\sqrt{19}},\ \sqrt{\tfrac32},\ \eta_6\Bigr\}.\]
Hence the set $\N$ has gaps between its smallest elements. The author also asked in~\cite[Section~1 and Question~4.9]{L14} whether $\N$ has non-empty interior, whether $[t,\infty)\subseteq\N$ for some $t$, and whether $\N$ is countable. In this paper, we answer these questions by establishing the following result.

\begin{theorem}\label{thm:main}
There is a real number $t$ such that $[t,\infty)\subseteq\N$. In particular, the set $\N$ has non-empty interior and is uncountable.
\end{theorem}

Theorem~\ref{thm:main} contrasts with the Toeplitz case. Let $B=[c_{i-j}]$ be a Toeplitz matrix indexed by $\bN\times\bN$. Then $B$ is a Schur multiplier if and only if there is a measure $\mu$ on the circle with $c_n=\hat\mu(n)$ for all $n\in\bZ$, and then $\sn B=\|\mu\|$, the total variation norm of $\mu$~\cite[Theorem~8.1]{Bennett}. If $B$ is also a $0$--$1$ matrix, then $\mu$ is idempotent under convolution, and the idempotent theorem of Helson and Cohen~\cite{Helson,Cohen} shows that $\{n\in\bZ:c_n=1\}$ lies in the coset ring of $\bZ$, the Boolean algebra of subsets of $\bZ$ generated by the cosets of subgroups. Green and Sanders~\cite{GS} proved a quantitative version of this theorem. The coset ring of $\bZ$ is countable, so the norms of the idempotent Toeplitz Schur multipliers form a countable set. Saeki~\cite{Saeki1,Saeki2} and Mudge and Pham~\cite{MP} determined the smallest elements of this set by classifying the idempotent measures of small norm.

The Schur norm of a matrix is also its factorisation norm $\gamma_2$ by a theorem that goes back to Grothendieck, see~\cite[Theorem~5.1]{Pisier}, and this norm is central in communication complexity. This has led to quantitative structural results for idempotent Schur multipliers, in the spirit of the theorem of Green and Sanders~\cite{HHH,GH,BHT}. Recently Beke, Goh, Hatami, Jaffe and Naylor~\cite{BGHJN} showed that every idempotent Schur multiplier is a finite signed sum of contractive idempotent Schur multipliers, answering a question of Katavolos and Paulsen~\cite[Section~2]{KatP}. Theorem~\ref{thm:main} shows that, unlike in the Toeplitz case, this structure theorem leaves room for a half-line of norms.

Idempotent Schur multipliers are also the idempotent bimodule maps over the diagonal masa of $\B(\ell^2)$, the maximal abelian self-adjoint subalgebra of diagonal operators~\cite[Section~1]{KatP}, and ranges of such bimodule projections are studied in the theory of reflexivity that goes back to Arveson~\cite{Arv,KatP,ET,ELT}. Parcet, de la Salle and Tablate~\cite{PST} characterise, locally and under a transversality condition, the idempotent Schur multipliers whose symbols are indicator functions of smooth domains and that are bounded on the Schatten $p$-classes for $1<p<\infty$ with $p\ne2$.

\subsection{The graphs}
\begin{figure}[ht]
\centering
\begin{tikzpicture}[x=1cm,y=1cm,font=\small,row/.style={circle,draw,fill=black,inner sep=1.5pt},col/.style={circle,draw,fill=white,inner sep=1.5pt}]
\draw[rounded corners,thick] (-1.5,-2.2) rectangle (1.5,2.2);
\node at (0,1.8) {$C$};
\foreach \i in {0,1,2} {
  \pgfmathsetmacro{\y}{1.5-1.5*\i}
  \foreach \s in {-1,1} {
    \pgfmathsetmacro{\yy}{\y+0.37*\s}
    \foreach \side in {-1,1} {
      \draw (\side*1,\y) -- (\side*2.4,\yy);
      \foreach \t in {-1,1} {
        \draw (\side*2.4,\yy) -- (\side*3.8,{\yy+0.17*\t});
      }
    }
  }
}
\foreach \i/\j in {0/1,0/2,1/0,1/2,2/0,2/1} {
  \pgfmathsetmacro{\yi}{1.5-1.5*\i}
  \pgfmathsetmacro{\yj}{1.5-1.5*\j}
  \draw[gray] (-1,\yi) -- (1,\yj);
}
\foreach \i in {0,1,2} {
  \pgfmathsetmacro{\y}{1.5-1.5*\i}
  \node[row] at (-1,\y) {};
  \node[col] at (1,\y) {};
  \foreach \s in {-1,1} {
    \pgfmathsetmacro{\yy}{\y+0.37*\s}
    \foreach \side/\first/\second in {-1/col/row,1/row/col} {
      \node[\first] at (\side*2.4,\yy) {};
      \foreach \t in {-1,1} {
        \node[\second] at (\side*3.8,{\yy+0.17*\t}) {};
      }
    }
  }
}
\foreach \x/\d in {-3.8/2,-2.4/1,-1/0,1/0,2.4/1,3.8/2} {
  \node at (\x,-2.6) {\scriptsize depth $\d$};
}
\end{tikzpicture}
\caption{The graph $T_C(m)$ for $C=J_3-I_3$ and $m=(2,2)$. Black vertices are rows and white vertices are columns.}
\label{fig:schematic}
\end{figure}
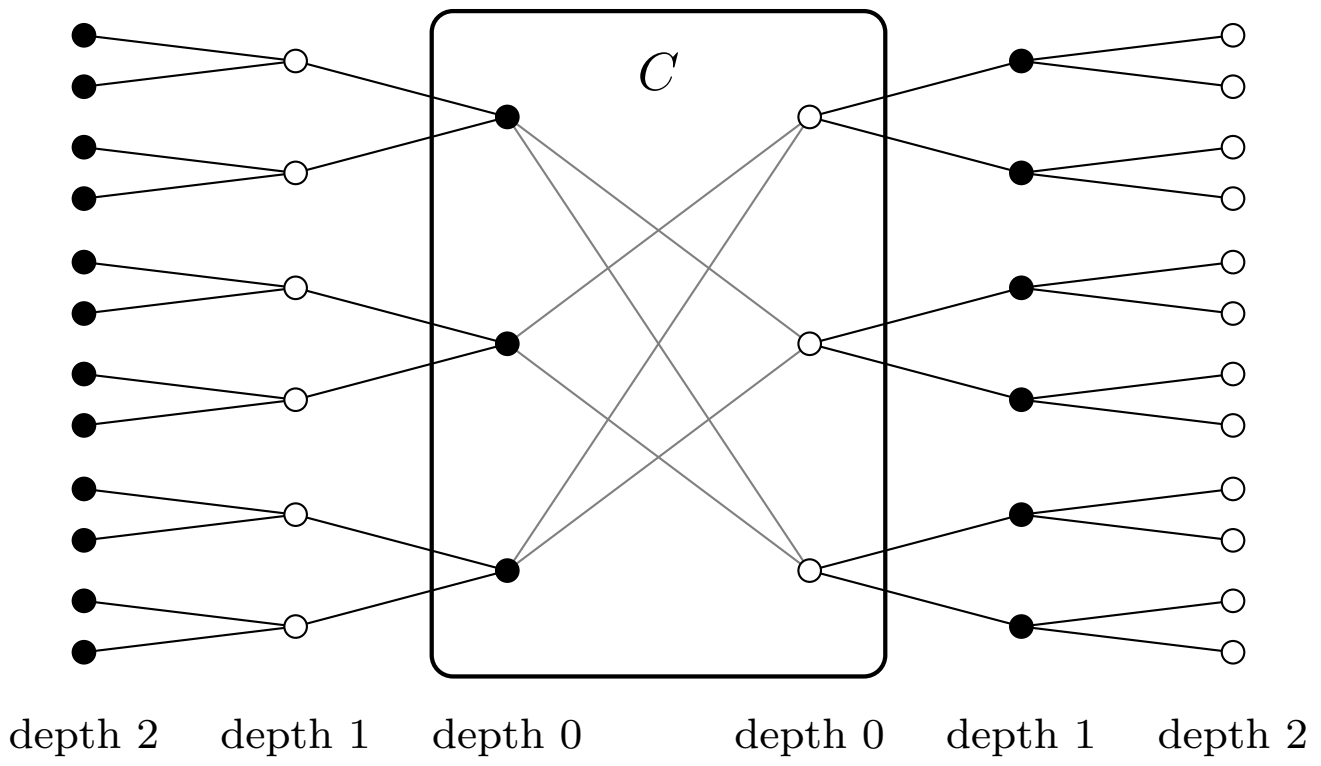

To prove Theorem~\ref{thm:main}, we realise a nontrivial interval as the norms of an explicit family of bipartite graphs, regarded as $0$--$1$ Schur multipliers. We build each graph from a certain finite bipartite graph $C$, which we call a \emph{core}, by hanging a rooted tree at every vertex of $C$. A core has $n$ row vertices and $n$ column vertices and a transitive group of automorphisms, see Definition~\ref{def:core}. The trees we hang on $C$ are specified by a finite or infinite sequence $m=(m_0,m_1,\dots)$ of positive integers, where $m_\ell$ is the number of children of every vertex at depth $\ell$. Here, the depth of a vertex is its distance from $C$ and the children of a vertex of depth $\ell$ are its neighbours of depth $\ell+1$. We call the $m_\ell$ the \emph{branching numbers} and write $T_C(m)$ for the resulting graph, as in Definition~\ref{def:trees} and Figure~\ref{fig:schematic}. A core $C$ is \emph{productive} if hanging sufficiently large trees in this way raises the norm strictly above both $2$ and $\sn C$, as made precise in Definition~\ref{def:productive}. Every productive core gives an inclusion $[t,\infty)\subseteq\N$ by Theorem~\ref{thm:halfline}. We identify a productive core in Section~\ref{sec:proof}, namely the bipartite graph
\begin{equation}\label{eq:core}C=(J_3-I_3)\otimes(J_4-I_4)\end{equation}
on $12+12$ vertices, where $I_k$ and $J_k$ are the $k\times k$ identity matrix and matrix of ones, respectively.

\subsection{Outline of the proof}
Theorem~\ref{thm:main} follows from three properties of the norms $\sn{T_C(m)}$, which for the present discussion we call~\ref{it:P1}, \ref{it:P2} and~\ref{it:P3}. The substance of the paper is to prove these properties for the graphs $T_C(m)$ built from a productive core $C$, and to find a productive core.

Throughout this discussion we take an integer $R$ that is large enough in terms of $C$, and a finite sequence $p$ with entries in $\Sigma_R=\{R,R+1,\dots,2R\}$, possibly the empty sequence $\varnothing$. Consider the closed intervals
\[\Lambda_p=\bigl[\sn{T_C(p,R,R,\dots)},\ \sn{T_C(p,2R,2R,\dots)}\bigr].\]
Since $\sn{T_C(m)}$ is non-decreasing in every $m_\ell$ by Lemma~\ref{lem:finite}, the norm $\sn{T_C(m)}$ lies in $\Lambda_p$ whenever $m$ begins with $p$ and has entries in $\Sigma_R$. In particular, we have $\Lambda_{(p,d)}\subseteq \Lambda_p$.

For a productive core $C$ and sufficiently large $R$, the three properties are the following.
\begin{enumerate}[label=(P\arabic*)]
\item\label{it:P1} The map $m\mapsto \sn{T_C(m)}$ is continuous on $\Sigma_R^\bN$ with respect to the product topology.
\item\label{it:P2} The interval $\Lambda_\varnothing$ has positive length.
\item\label{it:P3} The intervals $\Lambda_{(p,d)}$ and $\Lambda_{(p,d+1)}$ overlap, so that
\begin{equation}\label{eq:intro-union}
\Lambda_p=\bigcup_{d=R}^{2R}\Lambda_{(p,d)}.
\end{equation}
\end{enumerate}
These are established in Proposition~\ref{prop:tailchange} and parts~(a) and~(b) of Theorem~\ref{thm:large}, respectively.

We now give a sketch to see how Theorem~\ref{thm:main} follows. Let $y\in \Lambda_\varnothing$. By~\eqref{eq:intro-union}, we can successively choose $m_0,m_1,\dots$ with $y\in \Lambda_{(m_0,\dots,m_i)}$ for every $i$. Let $m=(m_0,m_1,\dots)$. We have $\sn{T_C(m)}\in \bigcap_{i\ge0} \Lambda_{(m_0,\dots,m_i)}$ and this intersection is a single point by~\ref{it:P1}, so $y=\sn{T_C(m)}$. Thus $\Lambda_\varnothing\subseteq\N$. Lemma~\ref{lem:interval} formalises this argument. The interval $\Lambda_\varnothing$ has positive length by~\ref{it:P2}, and the multiplicativity of Lemma~\ref{lem:tensor} turns it into $[t,\infty)$ for some $t$ by Lemma~\ref{lem:powers}. We remark that this digit-by-digit choice of $m$ resembles Shearer's solution~\cite{Shearer} of Hoffman's problem~\cite[Sections~1 and~6]{Hoffman} on limit points of spectral radii.

\subsection{Structure of the paper}
Section~\ref{sec:prelim} collects some general facts. The fundamental one is the dual formula of Lemma~\ref{lem:dual}, which states that for a finite matrix $B$, we have
\[\sn B=\max_w\,\|B\|_w\]
where $\|B\|_w$ is the trace norm of the matrix with the entries $\sqrt{w_iw_j}\,b_{ij}$, and $w$ runs over the functions from the rows and the columns to $[0,\infty)$ whose sums over the rows and over the columns are both $1$. We call the values of $w$ the \emph{weights} of the rows and columns, and the weights at which the maximum is attained \emph{optimal}. We also show by averaging that optimal weights may be assumed to be invariant under the automorphism group of a graph.

In Section~\ref{sec:trees}, we apply the dual formula to $T_C(m)$ to write the norm as the maximum of an expression involving trace norms of certain symmetric tridiagonal matrices of size at most the length of $m$, over weights that depend only on depth. In Section~\ref{sec:productive}, we define productive cores and show that optimal weights decay geometrically with depth, which yields the continuity~\ref{it:P1} in Proposition~\ref{prop:tailchange}. We also prove Theorem~\ref{thm:large}, which gives properties~\ref{it:P2} and~\ref{it:P3}.

Section~\ref{sec:proof} deduces the half-line of Theorem~\ref{thm:halfline}, using the covering Lemmas~\ref{lem:powers} and~\ref{lem:interval} of the outline. It also shows that the core of~\eqref{eq:core} is productive. In Section~\ref{sec:numerics}, we compute an explicit value of $t$, which is large. We also compare the constants of the proof with numerical computations and record open questions.

\section{Preliminaries}\label{sec:prelim}

The results of this section hold for general matrices and graphs, and several are known or follow quickly from known results. We include short arguments for the convenience of the reader. The following notation is used throughout the paper. For a vector $v$, we write $D_v$ for the diagonal matrix with diagonal $v$, and $\sqrt v$ is taken entrywise. For a finite sequence $p$ and a finite or infinite sequence $\nu$, we write $(p,\nu)$ for their concatenation, where an integer is identified with a sequence of length $1$, and for an integer $d$ we write $\overline{d}=(d,d,d,\dots)$. The trace norm of a matrix is $\|\cdot\|_1$, and for a finite vector $v$ we write $\|v\|_1=\sum_i|v_i|$ for its $\ell^1$ norm, which is the trace norm of $D_v$. The Hilbert--Schmidt norm is $\|\cdot\|_2$, the matrix unit $E_{ij}$ has entry $1$ at position $(i,j)$ and $0$ elsewhere, with its size clear from the context, and a matrix with countable index sets is identified with an operator between the corresponding $\ell^2$ spaces whenever it defines a bounded operator. A \emph{submatrix} of an $I\times J$ matrix $B$ is the restriction of $B$ to a product $I'\times J'$ with $I'\subseteq I$ and $J'\subseteq J$, and the corresponding graph is an \emph{induced subgraph} when $B$ is a $0$--$1$ matrix.

The following lemma is elementary, and its second assertion is~\cite[Proposition~3.1(3)]{L14}.
\begin{lemma}\label{lem:finite}
Let $B$ be a matrix indexed by countable sets $I\times J$. Then
\[\sn B=\sup\{\sn{B'}:B'\text{ is a finite submatrix of }B\}.\]
In particular, we have $\sn{B'}\le\sn B$ for every submatrix $B'$ of $B$.
\end{lemma}

We now show that the Schur norm of a finite matrix is a maximum over weights.

\begin{definition}\label{def:dual}
Let $B=[b_{ij}]$ be a finite matrix indexed by $I\times J$, where $I$ and $J$ are disjoint and non-empty. A \emph{weight} for $B$ is a function $w\colon I\sqcup J\to[0,\infty)$, and we call $w_i$ and $w_j$ the \emph{weights} of row $i$ and column $j$, respectively. For $S\subseteq I\sqcup J$ we write $w(S)=\sum_{v\in S}w_v$, and define
\[\Delta(I,J)=\{w\colon I\sqcup J\to[0,\infty)\;:\;w(I)=w(J)=1\}.\]
Given a weight $w$ for $B$ and a matrix $M=[M_{uv}]$ whose rows and columns are indexed by subsets of $I\sqcup J$, we define
\[\|M\|_w=\bigl\|\bigl[\sqrt{w_uw_v}\,M_{uv}\bigr]\bigr\|_1.\]
A \emph{maximiser} for $B$ is an element $w\in\Delta(I,J)$ such that
\[\|B\|_w\ge\|B\|_{w'}\qquad\text{for all }w'\in\Delta(I,J).\]
The weights of a maximiser are called \emph{optimal}. A maximiser for the matrix of a bipartite graph $G$ is also called a maximiser for $G$.
\end{definition}

\begin{lemma}\label{lem:dual}
Let $B$ be a finite matrix indexed by $I\times J$ for disjoint non-empty sets $I$ and $J$, and let $w$ be a weight for $B$.
\begin{enumerate}
\item[(a)] There is a maximiser for $B$, and
\[\sn B=\max\{\|B\|_{w'}:w'\in\Delta(I,J)\}.\]
\item[(b)] Let $A$ be the Hermitian matrix indexed by $(I\sqcup J)\times(I\sqcup J)$ given by
\[A=\begin{bmatrix}0&B\\B^*&0\end{bmatrix}.\]
Then $\|B\|_w=\tfrac12\|A\|_w$.
\item[(c)] We have $\|B\|_w\le\sqrt{w(I)w(J)}\,\sn B$.
\item[(d)] Let $B=B_1+B_2$ for matrices $B_1$ and $B_2$ indexed by $I\times J$. For $k=1,2$, let $I_k\subseteq I$ and $J_k\subseteq J$ be such that $B_k$ vanishes outside $I_k\times J_k$. Then
\[\|B\|_w\le\sqrt{w(I_1)w(J_1)}\,\sn{B_1}+\sqrt{w(I_2)w(J_2)}\,\sn{B_2}.\]
\end{enumerate}
\end{lemma}

\begin{proof}
(a) Under the bilinear pairing $\langle X,Y\rangle=\sum_{i,j}x_{ij}y_{ij}$, the dual of the operator norm is the trace norm, and the map $X\mapsto B\schur X$ is its own adjoint. So $\sn B$ is also the norm of $Y\mapsto B\schur Y$ on the trace class. Hence it is attained at an extreme point $\xi\eta^*$ of the unit ball of the trace class, where $\xi$ and $\eta$ are unit vectors. This is the argument of~\cite[Section~2]{Cowen}, where real matrices are treated.

Since $B\schur(\xi\eta^*)=D_\xi BD_{\overline{\eta}}$ and diagonal unitary matrices do not change the trace norm, we have
\[\|B\schur(\xi\eta^*)\|_1=\|B\|_{w'}\]
for the weight $w'\in\Delta(I,J)$ with $w'|_I=|\xi|^2$ and $w'|_J=|\eta|^2$, and every element of $\Delta(I,J)$ arises in this way.

(b) Let $W=D_{\sqrt{\smash[b]{w|_I}}\vphantom{\sqrt w_I}}BD_{\sqrt{\smash[b]{w|_J}}}$. Then $D_{\sqrt w}AD_{\sqrt w}=\left[\begin{smallmatrix}0&W\\W^*&0\end{smallmatrix}\right]$ and its square is $WW^*\oplus W^*W$. Since $WW^*$ and $W^*W$ have the same non-zero eigenvalues, we have $\|A\|_w=2\|W\|_1=2\|B\|_w$.

(c) If $w(I)=0$ or $w(J)=0$, then $w|_I=0$ or $w|_J=0$, and both sides vanish. Otherwise $\|B\|_w=\sqrt{w(I)w(J)}\,\|B\|_{w'}$ for the weight $w'\in\Delta(I,J)$ with $w'|_I=w|_I/w(I)$ and $w'|_J=w|_J/w(J)$, and (a) gives the claim.

(d) We have $\|B\|_w\le\|B_1\|_w+\|B_2\|_w$ by the triangle inequality. For $k=1,2$, let $w_k$ be the weight that agrees with $w$ on $I_k\sqcup J_k$ and vanishes elsewhere. Since $B_k$ vanishes outside $I_k\times J_k$, we have $\|B_k\|_w=\|B_k\|_{w_k}$, and $w_k$ has the totals $w(I_k)$ and $w(J_k)$ over the rows and over the columns. So (c) bounds the $k$th term by $\sqrt{w(I_k)w(J_k)}\,\sn{B_k}$, and adding these bounds establishes the claim.
\end{proof}

\begin{definition}\label{def:aut}
Let $G=(I,J,E)$ be a finite bipartite graph. An \emph{automorphism} of $G$ is an isomorphism of bipartite graphs from $G$ to itself in the sense of~\cite[Section~2]{L14}, that is, a permutation $\alpha$ of the vertex set $I\sqcup J$ that preserves adjacency and satisfies $\{\alpha(I),\alpha(J)\}=\{I,J\}$. The automorphisms form a group, which we denote by $\Aut(G)$.
\end{definition}

We now show that the maximisers for a graph form a convex set, and that averaging a maximiser over the automorphisms of the graph yields a maximiser whose weights are invariant.

\begin{lemma}\label{lem:invariant}
If $G=(I,J,E)$ is a finite bipartite graph with $I$ and $J$ non-empty, then there is a maximiser $w$ for $G$ such that $w\circ\alpha=w$ for every $\alpha\in\Aut(G)$.
\end{lemma}

\begin{proof}
Let $B$ be the matrix of $G$. We first show that the function $w\mapsto\|B\|_w$ is concave on $\Delta(I,J)$. The trace norm of a finite matrix $X$ satisfies
\[\|X\|_1=\tfrac12\min\bigl\{\|S\|_2^2+\|T\|_2^2:X=S^*T\bigr\},\]
where the minimum is over all finite matrices $S$ and $T$ with the same number of rows. Indeed, H\"older's inequality and the inequality of arithmetic and geometric means imply that $\|S^*T\|_1\le\|S\|_2\|T\|_2\le\frac12(\|S\|_2^2+\|T\|_2^2)$. A singular value decomposition $X=\sum_n\sigma_nu_nv_n^*$ gives equality for the matrices $S$ and $T$ with the rows $\sqrt{\sigma_n}u_n^*$ and $\sqrt{\sigma_n}v_n^*$, respectively.

Let $w$ be positive. The substitution $S=YD_{\sqrt{w|_I}}$ and $T=ZD_{\sqrt{w|_J}}$ shows that
\[\|B\|_w=\tfrac12\min\Bigl\{\sum_iw_i\|y_i\|^2+\sum_jw_j\|z_j\|^2:B=Y^*Z\Bigr\},\]
where the minimum is over all finite matrices $Y$ and $Z$ with the same number of rows, and $y_i$ and $z_j$ are the columns of $Y$ and $Z$, respectively. A minimum of affine functions is concave, so $w\mapsto\|B\|_w$ is concave on the positive elements of $\Delta(I,J)$, and by continuity on $\Delta(I,J)$. Hence the maximisers for $G$ form a convex set $\mathcal M$, which is non-empty by Lemma~\ref{lem:dual}(a).

Let $A$ be the symmetric $(I\sqcup J)\times (I\sqcup J)$ adjacency matrix of $G$ and let $\alpha\in\Aut(G)$. If $w\in\Delta(I,J)$, then $w\circ\alpha\in\Delta(I,J)$, since $\alpha$ permutes $\{I,J\}$. Since $A_{\alpha u,\alpha v}=A_{uv}$, the matrix $D_{\sqrt{w\circ\alpha}}AD_{\sqrt{w\circ\alpha}}$ is obtained from $D_{\sqrt w}AD_{\sqrt w}$ by permuting its rows and columns simultaneously by $\alpha$. Hence $\|A\|_{w\circ\alpha}=\|A\|_w$, and so $\|B\|_{w\circ\alpha}=\|B\|_w$ by Lemma~\ref{lem:dual}(b), and so $w\in\mathcal M$ implies $w\circ\alpha\in\mathcal M$. Therefore, if $w\in\mathcal M$, then the average
\[\overline w=\frac1{|\Aut(G)|}\sum_{\alpha\in\Aut(G)}w\circ\alpha\]
lies in $\mathcal M$ by convexity, and $\overline w\circ\beta=\overline w$ for every $\beta\in\Aut(G)$, since $\Aut(G)$ is a group.
\end{proof}

\section{Trees on a core}\label{sec:trees}

In this section, we define the graphs $T_C(m)$ and prove formulae for their norms that hold for every core. Section~\ref{sec:productive} then estimates these norms for productive cores.

\subsection{Chains}
\begin{definition}\label{def:chain}
For $L\ge1$ and $c=(c_0,\dots,c_{L-1})\in[0,\infty)^L$, let $Q(c)$ be the following symmetric tridiagonal matrix of order $L+1$.
\[
Q(c)=
\begin{bmatrix}
0 & c_0 & 0 & \cdots & 0 \\
c_0 & 0 & c_1 & \ddots & \vdots \\
0 & c_1 & 0 & \ddots & 0 \\
\vdots & \ddots & \ddots & \ddots & c_{L-1} \\
0 & \cdots & 0 & c_{L-1} & 0
\end{bmatrix}.
\]
We call $Q(c)$ a \emph{chain}.
\end{definition}

We will need the following facts about chains.

\begin{lemma}\label{lem:chain}
Let $L\ge1$.
\begin{enumerate}
\item[(a)] Let $Q$ be a real symmetric tridiagonal matrix of order $L+1$ whose off-diagonal entries $c_0,\dots,c_{L-1}$ are non-negative and whose diagonal is arbitrary. Then $\|Q\|_1$ is non-decreasing in each $c_e$.
\item[(b)] If $Q$ and $Q'$ are real symmetric tridiagonal matrices of order $L+1$ with diagonals $d$ and $d'$ and with off-diagonal entries $c_e$ and $c'_e$, then $\|Q-Q'\|_1\le\sum_i|d_i-d'_i|+2\sum_e|c_e-c'_e|$. In particular, we have $\|Q(c)\|_1\le2\sum_ec_e$.
\item[(c)] Let $c_1=0$ if $L=1$, and read $\|Q(c_2,\dots,c_{L-1})\|_1$ as $0$ if $L\le2$. Then
\[2c_0\le\|Q(c)\|_1-\|Q(c_2,\dots,c_{L-1})\|_1\le2\sqrt{c_0^2+c_1^2}.\]
\end{enumerate}
\end{lemma}

\begin{proof}
(a) Let $0\le\lambda<1$, and let $Q'$ be the matrix $Q$ with $c_e$ replaced by $\lambda c_e$. Let $D$ be the diagonal matrix with entries $1$ at the first $e+1$ positions and $-1$ at the remaining positions. Conjugation by $D$ changes the sign of the two entries $c_e$ and fixes all other entries of $Q$, so
\[Q'=\tfrac{1+\lambda}2\,Q+\tfrac{1-\lambda}2\,DQD.\]
Since $D$ is unitary, the triangle inequality gives $\|Q'\|_1\le\|Q\|_1$.

(b) We have $Q-Q'=\sum_i(d_i-d'_i)E_{ii}+\sum_e(c_e-c'_e)(E_{e+1,e+2}+E_{e+2,e+1})$. Each matrix $E_{ii}$ has trace norm $1$, and each matrix $E_{e+1,e+2}+E_{e+2,e+1}$ has trace norm $2$.

(c) Appending zeros to $c$ adds zero rows and columns to $Q(c)$, and it replaces a missing $Q(c_2,\dots,c_{L-1})$ by a zero matrix, which has trace norm $0$ as in the statement. So we may assume $L\ge3$. The pinching of $Q(c)$ to the first two coordinates and the last $L-1$ coordinates is
\[\begin{bmatrix}0&c_0\\c_0&0\end{bmatrix}\oplus Q(c_2,\dots,c_{L-1}),\]
and its trace norm is $2c_0+\|Q(c_2,\dots,c_{L-1})\|_1$. This is at most $\|Q(c)\|_1$ by the pinching inequality~\cite[(IV.52)]{Bhatia}, which is the lower bound. For the upper bound, we have
\[Q(c)=Q(c_0,c_1,0,\dots,0)+Q(0,0,c_2,\dots,c_{L-1}).\]
Since $\|Q(c_0,c_1)\|_1=2\sqrt{c_0^2+c_1^2}$, an application of the triangle inequality completes the upper bound.
\end{proof}

The next lemma reduces the trace norm of an operator built from partial isometries between consecutive summands to trace norms of chains. The decomposition that underlies this formula is standard, see Rojo and Medina~\cite{RM} for weighted graphs with a tree attached at every vertex and de la Harpe~\cite[Proposition~4.6]{dlH} for spherically symmetric trees. We include a proof, since neither source allows edge weights that vanish.

\begin{lemma}\label{lem:decomp}
Let $h\ge1$, let $K=K_0\oplus\dots\oplus K_h$ be an orthogonal direct sum of finite-dimensional Hilbert spaces, and let $k_j=\dim K_j$. For $0\le j<h$, let $\iota_j$ be a partial isometry on $K$ with initial space $K_j$ and final space contained in $K_{j+1}$, and let $c_j\ge0$. Let $T_0$ be a self-adjoint operator on $K_0$, regarded as an operator on $K$ that vanishes on $K_1\oplus\dots\oplus K_h$, and let
\[H=T_0+\sum_{j=0}^{h-1}c_j(\iota_j+\iota_j^*).\]
For $0\le\ell<h$, let $Q_\ell=Q(c_\ell,\dots,c_{h-1})$. If $\lambda_1,\dots,\lambda_{k_0}$ are the eigenvalues of $T_0$, repeated according to multiplicity, then
\[\|H\|_1=\sum_{i=1}^{k_0}\|\lambda_iE_{11}+Q_0\|_1+\sum_{\ell=1}^{h-1}(k_\ell-k_{\ell-1})\,\|Q_\ell\|_1.\]
\end{lemma}

\begin{proof}

We first choose an orthonormal basis of $K$. Let $K_{-1}=\{0\}$ and $\iota_{-1}=0$. Writing $k_{-1}=0$, let $d_\ell=k_\ell-k_{\ell-1}$. For $0\le\ell\le h$, let $f_{1,\ell},\dots,f_{d_\ell,\ell}$ be an orthonormal basis of $K_\ell\ominus\iota_{\ell-1}K_{\ell-1}$, where for $\ell=0$ we choose eigenvectors with $T_0f_{i,0}=\lambda_if_{i,0}$ for $1\le i\le k_0$. Define \[f_{i,\ell,s}=(\iota_{s-1}\dots\iota_\ell)f_{i,\ell}\in K_s\] for $0\le \ell\le s\le h$ and $1\le i\le d_\ell$.
This defines a collection of
\[\sum_{0\le\ell\le s\le h}d_\ell = \sum_{0\le s\le h}\sum_{0\le\ell\le s}(k_\ell-k_{\ell-1})=\sum_{0\le s\le h}k_s = \dim K\]
vectors, which are orthonormal since each $\iota_j$ is isometric on $K_j$ and $f_{i,\ell}\perp\iota_{\ell-1}K_{\ell-1}$ by construction. Hence they form an orthonormal basis of $K$.

We now compute $H$ in this basis. For $0\le\ell\le h$ and $1\le i\le d_\ell$, let
\[ W_{i,\ell}=\spn\{f_{i,\ell,s}:\ell\le s\le h\}.\]
Each of these subspaces is invariant for $T_0$ and each $\iota_j$, since for $0\le j<h$, we have
\[ T_0 f_{i,\ell,s}=\delta_{\ell,0}\delta_{s,0}\lambda_if_{i,0,0}\quad\text{and}\quad \iota_j f_{i,\ell,s}=\delta_{j,s}f_{i,\ell,s+1}.\]
Since $K$ is the orthogonal direct sum of the subspaces $W_{i,\ell}$, each one is reducing for $T_0$ and $\iota_j$.
These equations also imply that the matrix of $H|_{W_{i,\ell}}$ in the basis $f_{i,\ell,\ell},\dots,f_{i,\ell,h}$ is
\[ \lambda_i E_{11}+Q_0\quad\text{if }\ell=0,\qquad Q_\ell\quad\text{if }1\le \ell \le h-1,\qquad 0\quad\text{if }\ell=h.\]
So $H$ is the orthogonal direct sum of these restrictions. The trace norm of a direct sum is the sum of the trace norms, so the formula follows.
\end{proof}

\subsection{The graphs}
\begin{definition}\label{def:core}
A \emph{core} is a finite bipartite graph $C$ with at least one row vertex and at least one column vertex, such that $\Aut(C)$ acts transitively on the vertices of $C$.
\end{definition}

Since some automorphism of a core maps a row vertex to a column vertex, it exchanges the set of rows with the set of columns, so these sets have the same cardinality, which we write as $n$. So the vertex set $V(C)$ of a core has size $2n$.

For an integer $\ell\ge0$, we write $[\ell]=\{1,2,\dots,\ell\}$, so that $[0]=\emptyset$.

\begin{definition}\label{def:trees}
Let $C$ be a core, and let $m=(m_0,\dots,m_{h-1})$ be a finite sequence of positive integers of length $h\ge1$. For $0\le k\le h$, let
\[V_k=V(C)\times[m_0]\times[m_1]\times\dots\times[m_{k-1}],\]
where $V_0=V(C)$. We write $d(v)=k$ for the \emph{depth} of $v=(v_0,\dots,v_k)\in V_k$. If $k\ge1$, then we write $v^-=(v_0,\dots,v_{k-1})\in V_{k-1}$ for the \emph{parent} of $v$. The graph $T_C(m)$ has vertex set $V=V_0\cup\dots\cup V_h$ and edge set
\[E=E(C)\cup\{vv^-:d(v)\ge1\}.\]
A vertex $v$ is a row if and only if $v_0$ is a row and $d(v)$ is even, or $v_0$ is a column and $d(v)$ is odd. For $0\le k\le h$, let
\[P_k=n\,m_0m_1\cdots m_{k-1},\]
where the product $m_0m_1\cdots m_{k-1}$ is empty for $k=0$, so that $P_0=n$. Then $|V_k|=2P_k$, and $V_k$ contains $P_k$ rows and $P_k$ columns.
For an infinite sequence $m=(m_0,m_1,\dots)$ of positive integers, the graph $T_C(m)$ is the union of the graphs $T_C(m^{[h]})$ for $h\ge1$, where
\[m^{[h]}=(m_0,\dots,m_{h-1}).\]
In that case, $P_k$ is defined for every $k\ge0$ by the same formula.
\end{definition}

\begin{figure}[ht]
\centering
\begin{tikzpicture}[x=0.85cm,y=1cm,font=\small]
\node (a) at (0,2) {$(a)$};
\node (a1) at (-2,1) {$(a,1)$};
\node (a2) at (2,1) {$(a,2)$};
\node (a11) at (-3,0) {$(a,1,1)$};
\node (a12) at (-1,0) {$(a,1,2)$};
\node (a21) at (1,0) {$(a,2,1)$};
\node (a22) at (3,0) {$(a,2,2)$};
\foreach \p/\c in {a/a1,a/a2,a1/a11,a1/a12,a2/a21,a2/a22}
  \draw (\p) -- (\c);
\node[anchor=west] at (4.7,2) {$1$ vertex of depth $0$};
\node[anchor=west] at (4.7,1) {$m_0=2$ vertices of depth $1$};
\node[anchor=west] at (4.7,0) {$m_0m_1=4$ vertices of depth $2$};
\end{tikzpicture}
\caption{The vertices $v$ of $T_C(m)$ with $v_0=a$ for $m=(2,2)$. The parent of $(a,i,j)$ is $(a,i)$, and the parent of $(a,i)$ is $(a)$.}
\label{fig:tree}
\end{figure}
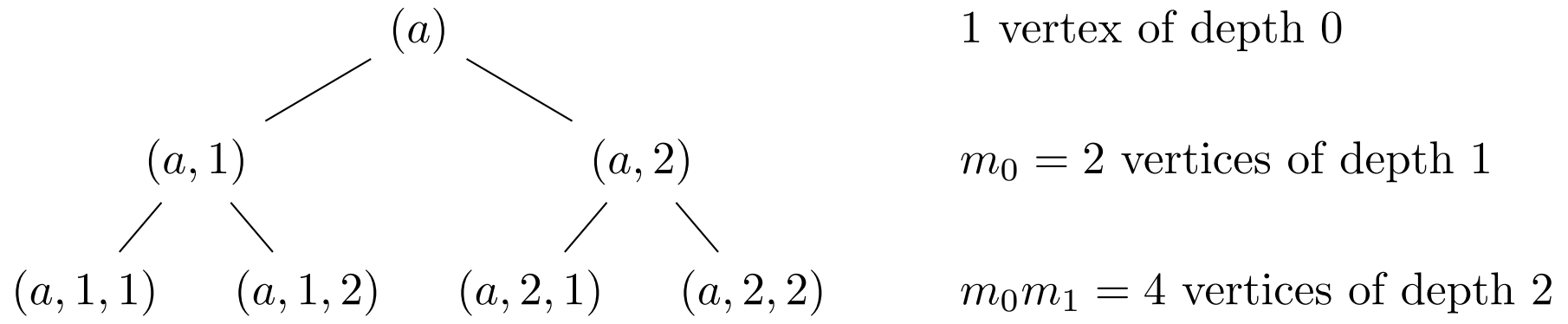

Thus $T_C(m)$ is obtained from $C$ by hanging at every vertex $a$ of $C$ the rooted tree on the vertices $v$ with $v_0=a$, in which every vertex of depth $\ell<h$ has $m_\ell$ children. Figure~\ref{fig:tree} shows the vertices $v$ with $v_0=a$ when $m=(2,2)$.

From now on until the end of Section~\ref{sec:proof}, we fix a core $C$ with $n$ row vertices and $n$ column vertices. Unless stated otherwise, $m=(m_0,\dots,m_{h-1})$ is a sequence of $h$ positive integers, for some integer $h\ge1$. We write $e_0,\dots,e_h$ for the standard basis of $\mathbb R^{h+1}$.

If $m$ has length at most that of $m'$, finite or infinite, and $m_\ell\le m'_\ell$ for every $\ell$ for which $m_\ell$ is defined, then $T_C(m)$ is an induced subgraph of $T_C(m')$. For an infinite sequence $m$, Lemma~\ref{lem:finite} shows that $\sn{T_C(m)}=\lim_{h\to\infty}\sn{T_C(m^{[h]})}$, since every finite subgraph of $T_C(m)$ lies in some $T_C(m^{[h]})$ and these norms increase with $h$.

\subsection{The norm formula}
In this subsection we prove Theorem~\ref{thm:norm}, which provides the formula for the norm of $T_C(m)$ used throughout the paper. It reduces the norm to a maximisation over $h+1$ weights, one for each depth, of an explicit function of the singular values of the core and of the trace norms of chains. The size of these chains is bounded by tree height rather than vertex count, so the formula makes it feasible to estimate the norms of highly branching trees on a core numerically. Moreover, it is key to our subsequent analysis. We deduce it from Lemma~\ref{lem:decomp}.

Let $B_C$ be the $n\times n$ matrix of $C$, and let $s_1,\dots,s_n$ be the singular values of $\tfrac1nB_C$, repeated according to multiplicity. We claim that
\begin{equation}\label{eq:snc}
s_1+\dots+s_n=\tfrac1n\|B_C\|_1=\sn C.
\end{equation}
Indeed, Lemma~\ref{lem:invariant} gives a maximiser for $C$ that is invariant under $\Aut(C)$, so its weights are constant on the vertices of $C$ because $\Aut(C)$ is transitive. Hence this maximiser is the weight $w$ that takes the value $\frac1n$ at every vertex of $C$, so that $\|B_C\|_w=\frac1n\|B_C\|_1$, and Lemma~\ref{lem:dual}(a) gives the claim. We remark that this is also the equality case of~\cite[Lemma~2.4]{Mathias}, which requires the diagonals of $|B_C|$ and $|B_C^*|$ to be constant. It is not difficult to see that this follows from the transitivity of $\Aut C$.

For $h\ge1$, let \[\Delta_h=\{x\in[0,\infty)^{h+1}:\|x\|_1=1\}\] be the $h$-dimensional probability simplex.

\begin{definition}\label{def:FM}
Write $T_C(m)=(I,J,E)$, and let $B$ be its matrix, indexed by $I\times J$. For $x=(x_0,\dots,x_h)\in[0,\infty)^{h+1}$, let $w$ be the weight for $B$ that takes the value $x_{d(v)}/P_{d(v)}$ at every vertex $v$ of $T_C(m)$. We call $w$ the \emph{weight induced by $x$}, and we define
\[F_m(x)=\|B\|_w.\]
A \emph{symmetric maximiser} for $m$ is a point $x\in\Delta_h$ at which $F_m$ attains its maximum on $\Delta_h$.
\end{definition}

By Definition~\ref{def:trees}, there are $P_k$ rows and $P_k$ columns at depth $k$, so the weight $w$ induced by $x$ gives the rows at depth $k$ total weight $x_k$, and similarly for the columns. Hence we have $w(I)=w(J)=\|x\|_1$, and if $x\in\Delta_h$, then $w\in \Delta(I,J)$. By Lemma~\ref{lem:dual}(c), it follows that
\begin{equation}\label{eq:fbound}
F_m(x)\le\|x\|_1\,\sn{T_C(m)}\qquad\text{for every }x\in[0,\infty)^{h+1}.
\end{equation}

\begin{definition}\label{def:rho}
Let $x=(x_0,\dots,x_h)\in[0,\infty)^{h+1}$. Let $c_k=\sqrt{x_kx_{k+1}}/P_k$ for $0\le k<h$, and let $Q_\ell=Q(c_\ell,\dots,c_{h-1})$ for $0\le\ell<h$. Let $n_0=n$ and $n_\ell=P_\ell-P_{\ell-1}$ for $1\le\ell<h$. Let
\[\rho_0=\frac1n\sum_{i=1}^n\|s_ix_0E_{11}+Q_0\|_1\qquad\text{and}\qquad \rho_\ell=\|Q_\ell\|_1\quad\text{for }1\le\ell<h.\]
We write $c_k(x)$, $Q_\ell(x)$ and $\rho_\ell(x)$ when the point $x$ needs disambiguation.
\end{definition}

By telescoping,
\begin{equation}\label{eq:fsum}
n_0+\dots+n_k=P_k\qquad\text{for }0\le k<h.
\end{equation}

\begin{theorem}\label{thm:norm}
Let $C$ be a core and let $m=(m_0,\dots,m_{h-1})$ be a sequence of positive integers.
\begin{enumerate}
\item[(a)] For every $x\in[0,\infty)^{h+1}$, we have
\[F_m(x)=\sum_{\ell=0}^{h-1}n_\ell\,\rho_\ell(x)=\sum_{i=1}^n\|s_ix_0E_{11}+Q_0(x)\|_1+\sum_{\ell=1}^{h-1}(P_\ell-P_{\ell-1})\,\|Q_\ell(x)\|_1.\]
\item[(b)] We have $\sn{T_C(m)}=\max_{x\in\Delta_h}F_m(x)$. If $x$ is a symmetric maximiser for $m$, then the weight induced by $x$ is a maximiser for $T_C(m)$.
\end{enumerate}
\end{theorem}

\begin{proof}
(a) We apply Lemma~\ref{lem:decomp} to the following matrix $H$. Let $V$ be the vertex set of $T_C(m)$, and let $E_{v,y}$ be the matrix units on $\mathbb{C}^V$. For $0\le k<h$, let
\[\iota_k=\frac1{\sqrt{m_k}}\sum_{v\in V_{k+1}}E_{v,v^-}.\]
Then $\iota_k$ is a partial isometry with initial space $\mathbb{C}^{V_k}$ and final space $\iota_k\mathbb{C}^{V_k}\subseteq\mathbb{C}^{V_{k+1}}$. Let $A_C$ be the adjacency matrix of $C$, acting on $\mathbb{C}^{V_0}$ and vanishing outside it. By Definition~\ref{def:trees}, the adjacency matrix of $T_C(m)$ is
\[A=A_C+\sum_{k=0}^{h-1}\sqrt{m_k}\,(\iota_k+\iota_k^*).\]
Let $w$ be the weight induced by $x$, which takes the value $x_k/P_k$ at the vertices of depth $k$. Let $D=D_{\sqrt w}$. Since $\iota_k$ maps the vertices of depth $k$ to those of depth $k+1$, we have $D\iota_kD=\sqrt{x_kx_{k+1}/(P_kP_{k+1})}\,\iota_k$, and $DA_CD=\frac{x_0}nA_C$. Since $P_{k+1}=m_kP_k$, we have
\[\sqrt{m_k}\sqrt{\frac{x_kx_{k+1}}{P_kP_{k+1}}}=\frac{\sqrt{x_kx_{k+1}}}{P_k}=c_k.\]
Hence the symmetric matrix
\[H=DAD=\frac{x_0}nA_C+\sum_{k=0}^{h-1}c_k\,(\iota_k+\iota_k^*)\]
satisfies $F_m(x)=\frac12\|A\|_w=\frac12\|H\|_1$ by Lemma~\ref{lem:dual}(b).

Let $K_k=\mathbb{C}^{V_k}$ for $0\le k\le h$, so that $\mathbb{C}^V=K_0\oplus\dots\oplus K_h$, and let $T_0=\tfrac1n\,x_0A_C$ on $K_0$. Then $H$ is the operator of Lemma~\ref{lem:decomp} for this data, with $k_j=|V_j|=2P_j$, so that $k_\ell-k_{\ell-1}=2n_\ell$ for $1\le\ell<h$. The adjacency matrix of $C$ is
\[A_C=\begin{bmatrix}0&B_C\\B_C^T&0\end{bmatrix},\]
so $T_0$ has the $2n$ eigenvalues $\pm x_0s_i$ for $1\le i\le n$, repeated according to multiplicity. Moreover, writing $S=\mathrm{diag}(1,-1,1,-1,\dots)$, we have $SE_{11}S=E_{11}$ and $SQ_0S=-Q_0$, so
\[\|-tE_{11}+Q_0\|_1=\|{-S}(tE_{11}+Q_0)S\|_1=\|tE_{11}+Q_0\|_1\qquad\text{for }t\in\mathbb R.\]
Hence Lemma~\ref{lem:decomp} implies that
\[F_m(x)=\tfrac12\|H\|_1=\sum_{i=1}^n\|s_ix_0E_{11}+Q_0\|_1+\sum_{\ell=1}^{h-1}n_\ell\|Q_\ell\|_1.\]
The first sum is $n\rho_0=n_0\rho_0$, so this is (a).

(b) By~\eqref{eq:fbound}, we have $F_m(x)\le\sn{T_C(m)}$ for every $x\in\Delta_h$. Conversely, let $V=V(T_C(m))$ and let $\Gamma=\Aut(T_C(m))$. For $\alpha\in\Aut(C)$ and permutations $\pi_i$ of $[m_{i-1}]$ for $1\le i\le h$, the map $\varphi:V\to V$ given by
\[\varphi(v_0,v_1,\dots,v_d)=(\alpha(v_0),\pi_1(v_1),\dots,\pi_d(v_d))\]
lies in $\Gamma$ because it commutes with taking parents and it preserves or exchanges the rows and the columns, and the row/column status of a vertex is determined by its first entry $v_0$ and the parity of its depth.  Since $C$ is a core, given $v,\hat v\in V_k$, we may choose $\alpha$ and $\pi_i$ with $\varphi(v)=\hat v$. So $\Gamma$ acts transitively on $V_k$ for every $k$.

By Lemma~\ref{lem:invariant}, there is a maximiser $w$ for $T_C(m)$ with $w\circ\gamma=w$ for every $\gamma\in\Gamma$. Hence $w$ is constant on $V_k$ for every $k$, by the transitivity just shown. Write this constant as $x_k/P_k$, so that $w$ is the weight induced by $x$. The rows of depth $k$ carry total weight $x_k$, so $x\in\Delta_h$, and by Lemma~\ref{lem:dual}(a) we have
\[F_m(x)=\|B\|_w=\sn{T_C(m)}.\]
This proves the first statement. For the second, let $x$ be a symmetric maximiser, and let $w$ be the weight induced by $x$. Then $w\in\Delta(I,J)$, and $\|B\|_w=F_m(x)=\sn{T_C(m)}$ by the first statement. So $w$ is a maximiser by Lemma~\ref{lem:dual}(a).
\end{proof}

\subsection{Bounds for \texorpdfstring{$F_m$}{F\_m}}
We will need the following Lipschitz bound for $F_m$.

\begin{lemma}\label{lem:lip}
Let $x,y\in[0,\infty)^{h+1}$. Then
\begin{equation}\label{eq:lip-main}|F_m(x)-F_m(y)|\le\sn C\,|x_0-y_0|+2\sum_{j=0}^{h-1}\bigl|\sqrt{x_jx_{j+1}}-\sqrt{y_jy_{j+1}}\bigr|.\end{equation}
\end{lemma}

\begin{proof}
By Lemma~\ref{lem:chain}(b), we have $|\rho_\ell(x)-\rho_\ell(y)|\le2\sum_{j\ge\ell}|c_j(x)-c_j(y)|$ for $\ell\ge1$. For each of the $n$ matrices $s_ix_0E_{11}+Q_0$, the same lemma bounds the difference of the trace norms at $x$ and $y$ by $s_i|x_0-y_0|+2\sum_j|c_j(x)-c_j(y)|$. Since $n_0=n$, summing over $i$ gives
\[n_0|\rho_0(x)-\rho_0(y)|\le(s_1+\dots+s_n)|x_0-y_0|+2n_0\sum_{j=0}^{h-1}|c_j(x)-c_j(y)|.\]
Adding the bounds for $\ell\ge1$, multiplied by $n_\ell$, Theorem~\ref{thm:norm}(a) gives
\[|F_m(x)-F_m(y)|\le(s_1+\dots+s_n)|x_0-y_0|+2\sum_{j=0}^{h-1}P_j|c_j(x)-c_j(y)|,\]
since $c_j$ enters $\rho_\ell$ only for $\ell\le j$, and $n_0+\dots+n_j=P_j$ by~\eqref{eq:fsum}. This is~\eqref{eq:lip-main}, since $s_1+\dots+s_n=\sn C$ by~\eqref{eq:snc}, and $P_jc_j(x)=\sqrt{x_jx_{j+1}}$.
\end{proof}

\section{Productive cores}\label{sec:productive}

In this section, we define productive cores and estimate the norm of $T_C(m)$ when $C$ is productive. Assuming productivity holds, we show that optimal weights decay geometrically with depth and that the norm is continuous for sufficiently large branching numbers, and we prove Theorem~\ref{thm:large}, which establishes positive interval length and the overlap of sub-intervals.

We will use the following bound on the Schur norm of two kinds of bipartite graph. A \emph{block} is a complete bipartite graph, and a \emph{direct sum of blocks} is a graph whose components are blocks or isolated vertices. A \emph{forest} is a graph without cycles. Every component of a forest is a tree, and we arbitrarily designate a vertex of every component as the root of that tree. This gives a parent to every vertex of the forest other than a root, namely its neighbour on the path to a root. The other neighbours of a vertex are its children.

We include a proof of the following elementary observation for completeness.
\begin{lemma}\label{lem:forest}
For every direct sum of blocks $M$ and every forest $F$ on countable vertex sets, we have $\sn M\le1$ and $\sn F\le2$.
\end{lemma}

\begin{proof}
The norm $\sn M$ is the supremum of the norms of the components of $M$, see~\cite[Proposition~3.2(3)]{L14}. These norms are $0$ or $1$, since the matrix of a block has all entries equal to $1$, so $\sn M\le1$.

Every edge of $F$ joins a vertex $v$ to its parent $v^-$. Let $F_0$ and $F_1$ be the graphs on the vertex set of $F$ whose edges $vv^-$ are those of $F$ for which the distance from the root of its component to $v$ is odd and even, respectively. This edge decomposition realises the matrix of $F$ as the sum of the matrices of $F_0$ and $F_1$. Every vertex has at most one parent, so every component of $F_0$ and of $F_1$ is an isolated vertex or a star, consisting of a vertex and its children. A star is complete bipartite, so $F_0$ and $F_1$ are direct sums of blocks. By the first part, $\sn F\le\sn{F_0}+\sn{F_1}\le1+1=2$.
\end{proof}

\subsection{Decay of the weights and continuity}
Proposition~\ref{prop:decay} below shows that the optimal weights decay geometrically with the depth, provided that the trees raise the norm of the core above $\sn C$ and above $2$. The following definition makes this condition precise.

\begin{definition}\label{def:productive}
The core $C$ is \emph{productive} with parameters $R_0$ and $h_0$ if $R_0$ and $h_0$ are positive integers and
\[\sn{T_C((\overline{R_0})^{[h_0]})}>\max\{2,\sn C\}.\]
\end{definition}

If $C$ is productive with these parameters, then
\begin{equation}\label{eq:productive}
\sn{T_C(m)}\ge\sn{T_C((\overline{R_0})^{[h_0]})}
\end{equation}
for every finite $m$ of length $h\ge h_0$ with all $m_\ell\ge R_0$,
since $T_C((\overline{R_0})^{[h_0]})$ is an induced subgraph of $T_C(m)$. Moreover, $T_C(m)$ is an edge-disjoint union of $C$ and a forest, so Lemma~\ref{lem:forest} and the triangle inequality imply that
\begin{equation}\label{eq:range}
\sn C\le\sn{T_C(m)}\le\sn C+2
\end{equation}
for every $m$, finite or infinite. %

To bound the decay rate of the optimal weights, Proposition~\ref{prop:decay} proves a recurrence of the form $\gamma u_k\le u_{k-1}+u_{k+1}$ for the square roots $u_k$ of the weights $x_k$. The next lemma shows that this yields a geometric bound.

\begin{lemma}\label{lem:recurrence}
Let $r>0$ and $\gamma\ge r+\tfrac1r$, let $h\ge1$, and let $u_0,\dots,u_h\ge0$ and $u_{h+1}=0$. Suppose that $\gamma u_k\le u_{k-1}+u_{k+1}$ for every $k$ with $1\le k\le h$ and $u_k>0$. Then $u_k\le ru_{k-1}$ for $1\le k\le h$.
\end{lemma}

\begin{proof}
Suppose that $u_j>ru_{j-1}$ for some $j$ with $1\le j\le h$. If $u_k>ru_{k-1}$ and $k<h$, then $u_k>0$, and so $u_{k+1}\ge\gamma u_k-u_{k-1}>(r+\tfrac1r)u_k-\tfrac1ru_k=ru_k$. So $u_h>ru_{h-1}$ by induction, and $u_h>0$. We have $\gamma u_h\le u_{h-1}<\tfrac1ru_h$ since $u_{h+1}=0$, and so $\gamma<\tfrac1r$, which contradicts $\gamma\ge r+\tfrac1r$. Hence $u_k\le ru_{k-1}$ for all $1\le k\le h$.
\end{proof}

\begin{proposition}\label{prop:decay}
Let $C$ be productive with parameters $R_0$ and $h_0$. Then there is $r_0\in(0,1)$ such that, for every $h\ge h_0$, every symmetric maximiser $x$ for a sequence $m_0,\dots,m_{h-1}$ of integers at least $R_0$ satisfies
\[x_k\le r_0^2x_{k-1}\quad\text{for }1\le k\le h.\]
In particular, $x_k\le r_0^{2k}$ for $0\le k\le h$.
\end{proposition}

\begin{proof}
Let $\gamma_0 =\sn{T_C((\overline{R_0})^{[h_0]})}$ and let $r_0\in(0,1)$ be the smaller root of the quadratic
\[r^2-\gamma_0 r+1=0,\]
which has real roots since $\gamma_0>2$. Let $\gamma=\sn{T_C(m)}$, so that $\gamma\ge\gamma_0=r_0+\tfrac1{r_0}$ by~\eqref{eq:productive}. By Theorem~\ref{thm:norm}(b), we have $\gamma=F_m(x)$. Let $u_j=\sqrt{x_j}$ for $0\le j\le h$ and $u_{h+1}=0$. By Lemma~\ref{lem:recurrence}, to establish the required decay $u_k\le r_0u_{k-1}$, it suffices to show that
\begin{equation}\label{eq:treedecay}
\gamma u_k\le u_{k-1}+u_{k+1}
\end{equation}
for every $k$ with $1\le k\le h$ and $x_k>0$. We prove this by comparing two bounds for the loss
\[L=F_m(x)-F_m(y)\qquad\text{where }y=x-te_k\text{ and }0<t<x_k,\]
incurred by removing $t$ from the weight at depth $k$.

Since $F_m(y)\le\gamma\|y\|_1=F_m(x)-\gamma t$ by~\eqref{eq:fbound}, we have
\begin{equation}\label{eq:loss-lower}
L\ge\gamma t.
\end{equation}
Let $\theta=u_k-\sqrt{x_k-t}=t(u_k+\sqrt{x_k-t})^{-1}$. Passing from $x$ to $y$ decreases $\sqrt{x_{k-1}x_k}$ by $u_{k-1}\theta$ and, if $k<h$, also decreases $\sqrt{x_kx_{k+1}}$ by $u_{k+1}\theta$. The values of $x_0$ and the other products $x_jx_{j+1}$ are unaltered. By Lemma~\ref{lem:lip},
\begin{equation}\label{eq:loss-upper}
L\le2(u_{k-1}+u_{k+1})\theta.
\end{equation}
Comparing~\eqref{eq:loss-lower} and~\eqref{eq:loss-upper}, dividing by $t$ and letting $t\to0$ yields~\eqref{eq:treedecay}. Lemma~\ref{lem:recurrence} now implies that $u_k\le r_0u_{k-1}$ for $1\le k\le h$, that is, $x_k\le r_0^2x_{k-1}$. Since $x_0\le\|x\|_1=1$, this gives $x_k\le r_0^{2k}$.
\end{proof}

We now use the decay of the weights to obtain the continuity property~\ref{it:P1} of the introduction.

\begin{proposition}\label{prop:tailchange}
If $C$ is productive with parameters $R_0$ and $h_0$, then the map $m\mapsto\sn{T_C(m)}$ is continuous on $\{R_0,R_0+1,\dots\}^{\bN}$ with the product topology.
\end{proposition}

\begin{proof}
Let $r_0$ be as in Proposition~\ref{prop:decay}, and let $\tau_k=\frac{2}{1-r_0^2}r_0^{2k}$. Since $\tau_k\to0$ as $k\to\infty$, it suffices to show that
\[\bigl|\sn{T_C(m)}-\sn{T_C(m')}\bigr|\le\tau_k\]
whenever $k\ge1$ and $m$ and $m'$ agree in the first $k$ places. We begin with finite sequences. Let $h\ge\max\{h_0,k\}$, and let $m,m'\in\{R_0,R_0+1,\dots\}^h$ agree in the first $k$ places. Let $m^+_\ell=\max(m_\ell,m'_\ell)$ and $m^-_\ell=\min(m_\ell,m'_\ell)$. Then $T_C(m^-)$ is an induced subgraph of $T_C(m)$ and of $T_C(m')$, which are induced subgraphs of $T_C(m^+)$. So both norms lie between $\sn{T_C(m^-)}$ and $\sn{T_C(m^+)}$, and it suffices to show that
\[\sn{T_C(m^+)}-\sn{T_C(m^-)}\le\tau_k.\]

Let $x$ be a symmetric maximiser for $m^+$, and let $w$ be the weight induced by $x$, which is a maximiser for $T_C(m^+)$ by Theorem~\ref{thm:norm}(b). Let $\sigma=\sum_{j\ge k}x_j$. Proposition~\ref{prop:decay} implies that
\[\sigma\le\sum_{j=k}^hr_0^{2j}\le\frac{r_0^{2k}}{1-r_0^2}.\] Let $B$ be the $I\times J$ matrix of $T_C(m^+)=(I,J,E)$. Let $I_1$ and $J_1$ be the sets of rows and of columns of depth at most $k$, and let $I_2$ and $J_2$ be those of depth at least $k$. Let $B_1$ be obtained from $B$ by replacing every entry outside $I_1\times J_1$ with $0$, and let $B_2=B-B_1$. Note that $B_1$ is the matrix of $T_C(m_0^-,\dots,m_{k-1}^-)$, which is an induced subgraph of $T_C(m^-)$, so $\sn{B_1}\le\sn{T_C(m^-)}$. Moreover, $B_2$ is the matrix of a forest, so $\sn{B_2}\le2$ by Lemma~\ref{lem:forest}, and the rows of depth $j$ and the columns of depth $j$ each carry total weight $x_j$ in $w$, so $w(I_2)=w(J_2)=\sum_{j\ge k}x_j=\sigma$.

Since $w(I_1)\le w(I)=1$ and $w(J_1)\le w(J)=1$, and $w$ is a maximiser, Lemma~\ref{lem:dual}(d) applied to $B=B_1+B_2$ gives
\[\sn{T_C(m^+)}=\|B\|_w\le\sn{T_C(m^-)}+2\sigma.\]
This proves the bound for finite sequences. The result for infinite sequences follows by taking the limit $h\to\infty$, since the truncations to length $h\ge\max\{h_0,k\}$ agree in the first $k$ places and their norms converge by Lemma~\ref{lem:finite}.
\end{proof}

\subsection{Weights and branching numbers}
This subsection proves Theorem~\ref{thm:large}. Its proof compares two norms at a single symmetric maximiser $x$. Raising a branching number by $1$ increases $F_m(x)$ only a little, whereas doubling a later branching number increases $F_m(x)$ by an amount proportional to $\sqrt{x_kx_{k+1}}$ for the next two depths $k$ and $k+1$. Lemma~\ref{lem:digit} makes both statements precise. The two changes involve weights at neighbouring depths, and Lemma~\ref{lem:weights} compares these, by showing that the weight at each depth from $2$ on is at least a fixed fraction of the weight at the previous depth. It also bounds the weight at depth $1$ from below, which yields the positive length~\ref{it:P2} in Theorem~\ref{thm:large}(a).

Let $C$ be productive with parameters $R_0$ and $h_0$. From this point on we use the following constants:
\begin{gather*}
  g=\tfrac12\bigl(\sn{T_C((\overline{R_0})^{[h_0]})}-\max\{2,\sn{C}\}\bigr),\qquad\beta=\sn{C}+3,\\
  R_1=\max\{R_0,\,2\beta^2+5\}\quad\text{and}\quad\delta_R=\Bigl(1-\frac4R\Bigr)\frac{g^2}{\beta^3R}\quad\text{for }R\ge R_1.
\end{gather*}
Note that $g>0$ since $C$ is productive, and $\delta_R>0$ since $R_1>4$.

\begin{lemma}\label{lem:weights}
Let $R\ge R_1$ and $h\ge h_0$, let $m_0,\dots,m_{h-1}\ge R$ be integers, let $x$ be a symmetric maximiser for $m$, and let $u_k=\sqrt{x_k}$ for $0\le k\le h$ and $u_{h+1}=0$. Then the following hold.
\begin{enumerate}
\item[(a)] We have $u_1\ge g$.
\item[(b)] We have $u_k\ge u_{k-1}/\beta$ for $2\le k\le h$.
\end{enumerate}
In particular, we have $x_k>0$ for $0\le k\le h$.
\end{lemma}

\begin{proof}
Let $\gamma=\sn{T_C(m)}$. By Theorem~\ref{thm:norm}(b), we have $\gamma=F_m(x)$.

(a) Since $F_m(0)=0$, Lemma~\ref{lem:lip} with $y=0$ gives
\[\gamma\le x_0\sn{C}+2\sum_{j=0}^{h-1}u_ju_{j+1}.\]
Since $2u_ju_{j+1}\le x_j+x_{j+1}$ and each $x_j$ with $j\ge1$ occurs in at most two of these terms, and since $u_0\le1$, we have
\begin{align*}
\gamma&\le x_0\sn{C}+2u_0u_1+2\sum_{j\ge1}x_j\le x_0\sn{C}+2(1-x_0)+2u_1\\
&\le\max\{2,\sn{C}\}+2u_1.
\end{align*}
Since $\gamma\ge\sn{T_C((\overline{R_0})^{[h_0]})}$ by~\eqref{eq:productive}, we get $u_1\ge g$.

(b) Let $2\le k\le h$. We may assume that $u_{k-1}>0$. Set $\alpha=\frac{u_k}{u_{k-1}}$ and $\varepsilon=\frac1R$. We add weight at depth $k$. Let $y$ be obtained from $x$ by replacing $x_k$ by \[y_k=(u_k+\varepsilon u_{k-1})^2=u_{k-1}^2(\alpha+\varepsilon)^2=x_k+\varepsilon u_{k-1}^2(2\alpha+\varepsilon).\] We bound $F_m(y)-F_m(x)$ from above and from below.

By~\eqref{eq:fbound}, we have $F_m(y)\le\gamma\|y\|_1$, and so
\begin{equation}\label{eq:gain-upper}
F_m(y)-F_m(x)\le\gamma\bigl(\|y\|_1-\|x\|_1\bigr)=\gamma\varepsilon u_{k-1}^2(2\alpha+\varepsilon).
\end{equation}

For the lower bound, passing from $x$ to $y$ increases $c_{k-1}$, does not decrease $c_k$ if $k<h$, and changes no other $c_j$. By Lemma~\ref{lem:chain}(a), no $\rho_\ell$ decreases, so by Theorem~\ref{thm:norm}(a), we have
\[F_m(y)-F_m(x)\ge n_{k-1}\bigl(\rho_{k-1}(y)-\rho_{k-1}(x)\bigr),\]
where $n_{k-1}=P_{k-1}-P_{k-2}\ge(1-\varepsilon)P_{k-1}$, since $P_{k-1}=m_{k-2}P_{k-2}$. The chains $Q_{k-1}(x)$ and $Q_{k-1}(y)$ differ only in their first two entries. Hence Lemma~\ref{lem:chain}(c) gives
\[\rho_{k-1}(y)-\rho_{k-1}(x)\ge2c_{k-1}(y)-2\sqrt{c_{k-1}(x)^2+c_k(x)^2},\]
where $c_k(x)$ is read as $0$ if $k=h$. Here $P_{k-1}c_{k-1}(y)=u_{k-1}^2(\alpha+\varepsilon)$ and $P_{k-1}c_{k-1}(x)=\alpha u_{k-1}^2$, and $P_{k-1}c_k(x)=u_ku_{k+1}/m_{k-1}$ is at most $\alpha\varepsilon u_{k-1}^2$, since $u_{k+1}\le u_{k-1}$ by Proposition~\ref{prop:decay}. As $\sqrt{1+\varepsilon^2}\le1+\tfrac12\varepsilon^2$, it follows that
\begin{align*}
P_{k-1}\bigl(\rho_{k-1}(y)-\rho_{k-1}(x)\bigr)&\ge2u_{k-1}^2(\alpha+\varepsilon)-2\alpha u_{k-1}^2\Bigl(1+\frac12\varepsilon^2\Bigr)\\
&=\varepsilon u_{k-1}^2(2-\alpha\varepsilon).
\end{align*}
Since $\rho_{k-1}(y)\ge\rho_{k-1}(x)$, it follows that
\begin{equation}\label{eq:gain-lower}
F_m(y)-F_m(x)\ge(1-\varepsilon)\varepsilon u_{k-1}^2(2-\alpha\varepsilon).
\end{equation}
Comparing~\eqref{eq:gain-upper} and~\eqref{eq:gain-lower} and dividing by $\varepsilon u_{k-1}^2$, we get
\[\gamma(2\alpha+\varepsilon)\ge(1-\varepsilon)(2-\alpha\varepsilon).\]
As $\gamma\le\beta-1$ by~\eqref{eq:range}, expanding gives $2(\beta-1)\alpha \ge 2 - (\beta+1)\varepsilon - \alpha\varepsilon(1-\varepsilon)$. Since $\alpha\le1$ by Proposition~\ref{prop:decay}, this implies that
\[2(\beta-1)\alpha\ge 2-(\beta+2)\varepsilon\ge 2-\frac2\beta,\]
where the last step uses $\varepsilon=\frac1R\le\frac2{\beta(\beta+2)}$ since $R\ge2\beta^2+5\ge\beta(\beta+2)/2$. Dividing by $2(\beta-1)$, we obtain $\alpha\ge 1/\beta$, i.e., $u_k\ge u_{k-1}/\beta$.

Finally, $u_0\ge u_1\ge g>0$ by Proposition~\ref{prop:decay} and (a), and so $u_k>0$ for $0\le k\le h$ by (b) and induction on $k$.
\end{proof}

The proof of Theorem~\ref{thm:large} compares graphs whose branching numbers differ at a single index. The following lemma estimates the resulting change of $F_m$ at any positive weights that decrease with the depth.

\begin{lemma}\label{lem:digit}
Let $R\ge3$ and $0\le\ell\le h-2$, and let $m$ and $m'$ be sequences of length $h$ with entries at least $R$ that agree except at the index $\ell$, where $m_\ell<m'_\ell$. Let $x\in(0,\infty)^{h+1}$, let $u_k=\sqrt{x_k}$ for $0\le k\le h$, and suppose that $u_k\le u_{k-1}$ for $1\le k\le h$. Then the following hold.
\begin{enumerate}
\item[(a)] We have $F_{m'}(x)-F_m(x)\le\dfrac{2R}{R-1}\,u_{\ell+1}u_{\ell+2}\Bigl(\dfrac1{m_\ell}-\dfrac1{m'_\ell}\Bigr)$.
\item[(b)] If $\ell\ge1$ and $m'_\ell=2m_\ell$, then $F_{m'}(x)-F_m(x)\ge\Bigl(1-\dfrac4R\Bigr)\dfrac{u_{\ell+1}u_{\ell+2}}{m_\ell}$.
\end{enumerate}
\end{lemma}

\begin{proof}
We make some preliminary observations. Let $\vartheta=m'_\ell/m_\ell$, let $c_j$, $n_k$ and $\rho_k$ be as in Definition~\ref{def:rho} for $m$ and $x$, and let $c'_j$, $Q'_k$, $n'_k$ and $\rho'_k$ be the corresponding objects for $m'$ and $x$. Since $P_{j+1}=m_jP_j$, the numbers $P_j$ for $m'$ are $P_j$ for $j\le\ell$ and $\vartheta P_j$ for $j>\ell$. Hence
\[c'_j=c_j\quad\text{for }j\le\ell\qquad\text{and}\qquad c'_j=\frac1\vartheta c_j\quad\text{for }j>\ell,\]
and
\[n'_k=n_k\ \text{for }k\le\ell,\qquad n'_{\ell+1}=\vartheta P_{\ell+1}-P_\ell,\qquad n'_k=\vartheta n_k\ \text{for }k\ge\ell+2.\]
For $k\ge\ell+1$, the chain $Q_k$ involves only $c_j$ with $j>\ell$, so $\rho'_k=\tfrac1\vartheta\rho_k$. Hence $n'_k\rho'_k=n_k\rho_k$ for $k\ge\ell+2$, and since $n_{\ell+1}=P_{\ell+1}-P_\ell$, we have
\[n'_{\ell+1}\rho'_{\ell+1}-n_{\ell+1}\rho_{\ell+1}=P_\ell\Bigl(1-\frac1\vartheta\Bigr)\rho_{\ell+1}.\]
By Theorem~\ref{thm:norm}(a), it follows that
\begin{equation}\label{eq:digit}
F_{m'}(x)-F_m(x)=P_\ell\Bigl(1-\frac1\vartheta\Bigr)\rho_{\ell+1}-\sum_{k=0}^{\ell}n_k\bigl(\rho_k-\rho'_k\bigr).
\end{equation}

The following will also be useful. Since $c_j=u_ju_{j+1}/P_j$ and $P_{j+1}=m_jP_j$, we have
\[\frac{c_{j+1}}{c_j}=\frac{u_{j+2}}{m_ju_j}\le\frac1R\qquad\text{for }0\le j<h-1,\]
because $u_{j+2}\le u_j$ and $m_j\ge R$. Summing a geometric series gives
\begin{equation}\label{eq:geomsum}
\sum_{j\ge i}c_j\le\frac R{R-1}c_i\qquad\text{for }0\le i<h.
\end{equation}
Moreover, since $P_{\ell+1}=m_\ell P_\ell$, we have
\begin{equation}\label{eq:pc}
P_\ell c_{\ell+1}=\frac{u_{\ell+1}u_{\ell+2}}{m_\ell}.
\end{equation}

(a) Since $c'_j\le c_j$ for every $j$, Lemma~\ref{lem:chain}(a) shows that $\rho'_k\le\rho_k$ for $0\le k\le h$, so the sum in~\eqref{eq:digit} is non-negative. By Lemma~\ref{lem:chain}(b) and~\eqref{eq:geomsum}, we have
\[\rho_{\ell+1}\le2\sum_{j\ge\ell+1}c_j\le\frac{2R}{R-1}c_{\ell+1}.\]
Since $1-\tfrac1\vartheta=1-\tfrac{m_\ell}{m'_\ell}$, it follows from~\eqref{eq:digit} and~\eqref{eq:pc} that
\[F_{m'}(x)-F_m(x)\le\frac{2R}{R-1}P_\ell c_{\ell+1}\Bigl(1-\frac{m_\ell}{m'_\ell}\Bigr)=\frac{2R}{R-1}u_{\ell+1}u_{\ell+2}\Bigl(\frac1{m_\ell}-\frac1{m'_\ell}\Bigr).\]

(b) Now $\vartheta=2$. By Lemma~\ref{lem:chain}(c), we have $\rho_{\ell+1}\ge2c_{\ell+1}$, so the first term in~\eqref{eq:digit} is at least $P_\ell c_{\ell+1}$. We show below that
\begin{equation}\label{eq:digit-low}
\sum_{k=0}^{\ell-1}n_k\bigl(\rho_k-\rho'_k\bigr)\le\frac1{R-1}P_\ell c_{\ell+1}
\end{equation}
and
\begin{equation}\label{eq:digit-ell}
n_\ell\bigl(\rho_\ell-\rho'_\ell\bigr)\le\Bigl(\frac1R+\frac1{R-1}\Bigr)P_\ell c_{\ell+1}.
\end{equation}
Given these bounds, \eqref{eq:digit} implies that
\[F_{m'}(x)-F_m(x)\ge\Bigl(1-\frac1R-\frac2{R-1}\Bigr)P_\ell c_{\ell+1}\ge\Bigl(1-\frac4R\Bigr)P_\ell c_{\ell+1},\]
where the last step uses $2R\le3(R-1)$, which holds since $R\ge3$. By~\eqref{eq:pc}, this is (b).

We now prove~\eqref{eq:digit-low}. Let $k<\ell$. The entries $c_j$ and $c'_j$ of $Q_k$ and $Q'_k$ differ only for $j\ge\ell+1$, where $c_j-c'_j=\tfrac12c_j$. By the triangle inequality, Lemma~\ref{lem:chain}(b) and~\eqref{eq:geomsum}, we have
\[\rho_k-\rho'_k\le\sum_{j\ge\ell+1}c_j\le\frac R{R-1}c_{\ell+1},\]
where for $k=0$ we apply Lemma~\ref{lem:chain}(b) to each of the $n$ pairs $s_ix_0E_{11}+Q_0$ and $s_ix_0E_{11}+Q'_0$, which have the same diagonal, and average. By~\eqref{eq:fsum}, we have $n_0+\dots+n_{\ell-1}=P_{\ell-1}\le\tfrac1RP_\ell$, and~\eqref{eq:digit-low} follows.

Finally, we prove~\eqref{eq:digit-ell}. Since $\ell\ge1$, we have $\rho_\ell=\|Q_\ell\|_1$ and $\rho'_\ell=\|Q'_\ell\|_1$. Let $\tilde Q=Q(c_{\ell+2},\dots,c_{h-1})$, read as $0$ if $h=\ell+2$. Then
\[Q_\ell=Q(c_\ell,c_{\ell+1},\dots,c_{h-1})\qquad\text{and}\qquad Q'_\ell=Q\bigl(c_\ell,\tfrac12c_{\ell+1},\dots,\tfrac12c_{h-1}\bigr).\]
The upper bound of Lemma~\ref{lem:chain}(c) for $Q_\ell$ and its lower bound for $Q'_\ell$ give
\[\rho_\ell\le2\sqrt{c_\ell^2+c_{\ell+1}^2}+\|\tilde Q\|_1\qquad\text{and}\qquad\rho'_\ell\ge2c_\ell+\tfrac12\|\tilde Q\|_1.\]
Subtracting these bounds yields $\rho_\ell-\rho'_\ell\le(2\sqrt{\smash[b]{c_\ell^2+c_{\ell+1}^2}}-2c_\ell)+\tfrac12\|\tilde Q\|_1$.
Since $c_\ell>0$, we have
\[2\sqrt{c_\ell^2+c_{\ell+1}^2}-2c_\ell=\frac{2c_{\ell+1}^2}{\sqrt{c_\ell^2+c_{\ell+1}^2}+c_\ell}\le\frac{c_{\ell+1}}{c_\ell}\,c_{\ell+1}\le\frac1Rc_{\ell+1},\]
and by Lemma~\ref{lem:chain}(b) and~\eqref{eq:geomsum},
\[\tfrac12\|\tilde Q\|_1\le\sum_{j\ge\ell+2}c_j\le\Bigl(\frac R{R-1}-1\Bigr)c_{\ell+1}=\frac1{R-1}c_{\ell+1}.\]
Since $n_\ell\le P_\ell$, \eqref{eq:digit-ell} follows.
\end{proof}

We now use these lemmas to obtain the positive length property~\ref{it:P2} and the overlap property~\ref{it:P3} of the introduction, in parts (a) and (b) of the following theorem.

\begin{theorem}\label{thm:large}
Let $C$ be productive with parameters $R_0$ and $h_0$, and let $R\ge R_1$. Then the following hold, where all sequences have entries in $\{R,R+1,\dots,2R\}$.
\begin{enumerate}
\item[(a)] $\sn{T_C(\overline{2R})}-\sn{T_C(\overline{R})}\ge\delta_R>0$.
\item[(b)] For every finite sequence $p$ of length $\ell\ge0$ and every $d$ with $R\le d<2R$, we have $\sn{T_C(p,d+1,\overline{R})}\le\sn{T_C(p,d,\overline{2R})}$.
\end{enumerate}
\end{theorem}

\begin{proof}
In both parts, we take $h\ge h_0$ and all sequences are truncated to length $h$, so that the corresponding graphs are finite. Let $x$ be a symmetric maximiser for one of the truncated sequences, and let $u_k=\sqrt{x_k}$. By Lemma~\ref{lem:weights} and Proposition~\ref{prop:decay}, we have $x_k>0$ for $0\le k\le h$ and $u_k\le u_{k-1}$ for $1\le k\le h$, so Lemma~\ref{lem:digit} applies at $x$. By Theorem~\ref{thm:norm}(b), we have $\sn{T_C(m)}\ge F_m(x)$ for every truncated sequence $m$, with equality if $x$ is a symmetric maximiser for $m$.

(a) Let $h\ge\max\{h_0,3\}$, and let $x$ be a symmetric maximiser for $\overline R$. By Lemma~\ref{lem:weights}, we have $\beta^3u_2u_3\ge u_1^2\ge g^2$. So Lemma~\ref{lem:digit}(b) at the index $1$ gives
\[\sn{T_C(R,2R,\overline{R})}-\sn{T_C(\overline{R})}\ge F_{(R,2R,\overline{R})}(x)-F_{\overline R}(x)\ge\Bigl(1-\frac4R\Bigr)\frac{u_2u_3}R\ge\delta_R.\]
Letting $h\to\infty$ yields the same inequality for the infinite sequences. Since $T_C(R,2R,\overline{R})$ is an induced subgraph of $T_C(\overline{2R})$, this proves (a).

(b) Let $h\ge\max\{h_0,\ell+3\}$, and let $m=(p,d,\overline{R})$, $m'=(p,d+1,\overline{R})$ and $m''=(p,d,2R,\overline{R})$. Let $x$ be a symmetric maximiser for $m'$. Then
\[\sn{T_C(m'')}-\sn{T_C(m')}\ge\bigl(F_{m''}-F_m\bigr)(x)-\bigl(F_{m'}-F_m\bigr)(x).\]
We now bound the two differences using Lemma~\ref{lem:digit} at the indices $\ell+1$ and $\ell$. Note that $1\le\ell+1\le h-2$ since $h\ge\ell+3$. By Lemma~\ref{lem:digit}(b), the first difference is at least $(1-\tfrac4R)\tfrac1Ru_{\ell+2}u_{\ell+3}$. By Lemma~\ref{lem:weights}(b) with $k=\ell+2$ and $k=\ell+3$, we have $u_{\ell+1}u_{\ell+2}\le\beta^2u_{\ell+2}u_{\ell+3}$, and $\tfrac1d-\tfrac1{d+1}\le\tfrac1{R^2}$. By Lemma~\ref{lem:digit}(a), the second is at most $\tfrac{2\beta^2}{R-1}\tfrac1Ru_{\ell+2}u_{\ell+3}$. Hence
\[\sn{T_C(m'')}-\sn{T_C(m')}\ge\Bigl(1-\frac4R-\frac{2\beta^2}{R-1}\Bigr)\frac{u_{\ell+2}u_{\ell+3}}R.\]
We can conclude for finite sequences by observing that this lower bound is non-negative, because $R\ge2\beta^2+5$ and so $(1-\tfrac4R)(R-1)=R-5+\tfrac4R\ge2\beta^2$. Letting $h\to\infty$ gives the same inequality for the infinite sequences. Since $T_C(p,d,2R,\overline{R})$ is an induced subgraph of $T_C(p,d,\overline{2R})$, this proves (b).
\end{proof}

\section{Proof of the main theorem}\label{sec:proof}

Theorem~\ref{thm:halfline} below shows that the norms $\sn{T_C(m)}$ contain a half-line when $C$ is a productive core. It follows from Proposition~\ref{prop:tailchange}, from Theorem~\ref{thm:large}, which establishes positive length and the absence of gaps, and from the covering lemmas below. We then show that a specific core is productive, which proves Theorem~\ref{thm:main}.

The following lemmas contain the covering argument of the outline. The closure of $\N$ under multiplication in the first lemma is attributed to Davidson in~\cite[Remark~5]{KatP}.
\begin{lemma}\label{lem:tensor}
If $B$ and $B'$ are Schur multipliers indexed by countable sets, then their Kronecker product $B\otimes B'$ satisfies $\sn{B\otimes B'}=\sn B\,\sn{B'}$. Hence, the set $\N$ is closed under multiplication.
\end{lemma}

\begin{proof}
The inequality $\le$ is~\cite[Theorem~10.5]{Bennett}, and $\ge$ follows from $(B\otimes B')\schur(X\otimes X')=(B\schur X)\otimes(B'\schur X')$ since $\|X\otimes X'\|=\|X\|\,\|X'\|$. A bijection $\bN\times\bN\to\bN$ gives the last statement.
\end{proof}

\begin{lemma}\label{lem:powers}
Let $S\subseteq\bR$ be closed under multiplication, and let $1<a<b$ with $[a,b]\subseteq S$. If $K$ is an integer with $K\log(b/a)\ge\log a$, then $[a^K,\infty)\subseteq S$.
\end{lemma}

\begin{proof}
As $\log a>0$, we have $K\ge1$. Hence $[a^K,\infty)=\bigcup_{k\ge K}[a^k,b^k]\subseteq S$, since every point of $[a^k,b^k]$ is $s^k$ for some $s\in[a,b]$ and $S$ is closed under multiplication, and $b^k\ge a^{k+1}$ for $k\ge K$.
\end{proof}

\begin{lemma}\label{lem:interval}
Let $d_0\le d_1$ be integers, and let $\Sigma=\{d_0,d_0+1,\dots,d_1\}^{\bN}$ with the product topology. Let $\Phi\colon\Sigma\to\bR$ be continuous with $\Phi(\mu)\le\Phi(\nu)$ whenever $\mu_i\le\nu_i$ for every $i$. Suppose that for every finite sequence $p$ with entries in $\{d_0,d_0+1,\dots,d_1\}$, including the empty sequence $\varnothing$, and every $d$ with $d_0\le d<d_1$, we have
\begin{equation}\label{eq:interval-hyp}
\Phi(p,d+1,\overline{d_0})\le\Phi(p,d,\overline{d_1}).
\end{equation}
Then $\Phi(\Sigma)=\bigl[\Phi(\overline{d_0}),\Phi(\overline{d_1})\bigr]$.
\end{lemma}

\begin{proof}
For a finite sequence $p$ with entries in $\{d_0,\dots,d_1\}$, let
\[\Lambda_p=\bigl[\Phi(p,\overline{d_0}),\ \Phi(p,\overline{d_1})\bigr],\]
so that $\Lambda_\varnothing=[\Phi(\overline{d_0}),\Phi(\overline{d_1})]$. By monotonicity, $\Lambda_{(p,d)}\subseteq \Lambda_p$. The interval $\Lambda_{(p,d_0)}$ has the same left endpoint as $\Lambda_p$, the interval $\Lambda_{(p,d_1)}$ has the same right endpoint, and by~\eqref{eq:interval-hyp}, the interval $\Lambda_{(p,d+1)}$ starts no later than $\Lambda_{(p,d)}$ ends. Hence
\begin{equation}\label{eq:interval-union}
\Lambda_p=\bigcup_{d=d_0}^{d_1}\Lambda_{(p,d)}.
\end{equation}

Now let $y\in \Lambda_\varnothing$. By~\eqref{eq:interval-union}, we can choose $m_0,m_1,\dots$ one at a time so that $y\in \Lambda_{(m_0,\dots,m_i)}$ for every $i$. Let $m=(m_0,m_1,\dots)\in\Sigma$. The endpoints $\Phi(m_0,\dots,m_i,\overline{d_0})$ and $\Phi(m_0,\dots,m_i,\overline{d_1})$ of $\Lambda_{(m_0,\dots,m_i)}$ converge to $\Phi(m)$ as $i\to\infty$ by continuity, so $y=\Phi(m)$. This shows that $\Lambda_\varnothing\subseteq\Phi(\Sigma)$, and the reverse inclusion holds by monotonicity.
\end{proof}

\begin{theorem}\label{thm:halfline}
If $C$ is a productive core with parameters $R_0$ and $h_0$, then there is a real number $t$ with $[t,\infty)\subseteq\N$. In fact, for every integer $R\ge R_1$ one can take $t=(\sn{C}+2)^{K}$ with
\[K=\left\lceil\left(\frac{\log(\sn{C}+2+\delta_R)}{\log(\sn{C}+2)}-1\right)^{-1}\right\rceil,\]
where $R_1$ and $\delta_R$ are the constants of Section~\ref{sec:productive}.
\end{theorem}

\begin{proof}
Let $R\ge R_1$, let $\Sigma=\{R,R+1,\dots,2R\}^{\bN}$, and let $\Phi\colon\Sigma\to\bR$ be given by $\Phi(m)=\sn{T_C(m)}$.

By~\eqref{eq:range}, every graph $T_C(m)$ with $m\in\Sigma$ is a Schur multiplier of norm at most $\sn{C}+2$, so $\Phi$ takes values in $\N$. We now show that $\Phi$ satisfies the hypotheses of Lemma~\ref{lem:interval}. If $\mu,\nu\in\Sigma$ with $\mu_i\le\nu_i$ for every $i$, then $T_C(\mu)$ is an induced subgraph of $T_C(\nu)$, so Lemma~\ref{lem:finite} shows that $\Phi(\mu)\le\Phi(\nu)$. Since $R\ge R_1\ge R_0$, Proposition~\ref{prop:tailchange} shows that $\Phi$ is continuous. Finally, the conclusion of Theorem~\ref{thm:large}(b) is the hypothesis~\eqref{eq:interval-hyp}. By Lemma~\ref{lem:interval}, we have
\[[a,b]=\Phi(\Sigma)\subseteq\N,\qquad\text{where}\quad a=\sn{T_C(\overline{R})}\quad\text{and}\quad b=\sn{T_C(\overline{2R})}.\]

We now bound $a$ and $b$. Theorem~\ref{thm:large}(a) gives $b-a\ge\delta_R>0$, and $b\le\sn{C}+2$ by~\eqref{eq:range}. Applying~\eqref{eq:productive} to the sequences $\overline R^{[h]}$ with $h\ge h_0$ and using Lemma~\ref{lem:finite}, we obtain $a\ge\sn{T_C((\overline{R_0})^{[h_0]})}>2$. Hence $1<a<b\le\kappa$, where $\kappa=\sn{C}+2$. Since $b\ge a+\delta_R$ and $a\le\kappa$, we have
\[\frac ba\ge1+\frac{\delta_R}a\ge1+\frac{\delta_R}\kappa=\frac{\kappa+\delta_R}\kappa.\]
So $\log(b/a)\ge\log(\kappa+\delta_R)-\log\kappa$ and $\log a\le\log\kappa$, and the choice of $K$ implies that $K\log(b/a)\ge\log a$. The set $\N$ is closed under multiplication by Lemma~\ref{lem:tensor}, so Lemma~\ref{lem:powers} yields $[t,\infty)\subseteq[a^{K},\infty)\subseteq\N$.
\end{proof}

We will now show that $C=(J_3-I_3)\otimes(J_4-I_4)$ is a productive core, which completes the proof of Theorem~\ref{thm:main}.

\begin{proposition}\label{prop:cert}
The graph $C=(J_3-I_3)\otimes(J_4-I_4)$ is a core with $\sn C=2$, and it is productive with parameters $R_0=40$ and $h_0=6$.
\end{proposition}

\begin{proof}
The rows and the columns of $C$ are both indexed by $\{1,2,3\}\times\{1,2,3,4\}$, and $C$ is preserved by applying the same map $\sigma\times\pi$ to the rows and to the columns, for permutations $\sigma$ of $\{1,2,3\}$ and $\pi$ of $\{1,2,3,4\}$. Since the matrix of $C$ is symmetric, the graph $C$ is also preserved by the exchange of each row with the column of the same index. These automorphisms act transitively on the vertices, so $C$ is a core.

The matrices $J_3-I_3$ and $J_4-I_4$ have the eigenvalues $2,-1,-1$ and $3,-1,-1,-1$, respectively. So the singular values of $B_C$ are the absolute values of the twelve products of these eigenvalues, and $n=12$. Hence the singular values of $\tfrac1nB_C$ are
\[s_i=\frac{\mu_i}{12}\quad\text{where}\quad(\mu_1,\dots,\mu_{12})=(6,\;3,3,\;2,2,2,\;1,1,1,1,1,1).\]
The $\mu_i$ sum to $24$, so $\sn C=2$ by~\eqref{eq:snc}. Alternatively, $\sn C=\sn{J_3-I_3}\,\sn{J_4-I_4}=\tfrac43\cdot\tfrac32$ by~\cite[Theorem~2.1]{BCD} and Lemma~\ref{lem:tensor}.

It remains to show that $\sn{T_C(m)}>2.003$ for $m=(40,40,40,40,40,40)$. Let
\[u=(24,6,5,4,3,2,1)\qquad\text{and}\qquad x_k=u_k^2\quad\text{for }0\le k\le6,\]
so that $x=(576,36,25,16,9,4,1)$ and $\|x\|_1=667$. By~\eqref{eq:fbound}, it suffices to show that
\[\frac{F_m(x)}{\|x\|_1}>2.003.\]
By Theorem~\ref{thm:norm}(a), we have
\[F_m(x)=\sum_{\ell=0}^5n_\ell\rho_\ell(x),\qquad n_0=12\quad\text{and}\quad n_\ell=12\cdot39\cdot40^{\ell-1}\text{ for }1\le\ell\le5,\]
where $n_\ell$ and $\rho_\ell$ are as in Definition~\ref{def:rho}. At $x$ we have $c_k=u_ku_{k+1}/P_k$ with $P_k=12\cdot40^k$.

For $1\le\ell\le5$, Lemma~\ref{lem:chain}(c) gives $\rho_\ell(x)=\|Q_\ell\|_1\ge2c_\ell$, and $n_\ell c_\ell=\tfrac{39}{40}u_\ell u_{\ell+1}$. Since $\sum_{\ell=1}^5u_\ell u_{\ell+1}=70$, it follows that
\[\sum_{\ell=1}^5n_\ell\rho_\ell(x)\ge2\cdot\frac{39}{40}\cdot70=136.5.\]

For $\ell=0$, we have $c_0=12$, and $s_ix_0E_{11}+Q_0$ has diagonal entry $48\mu_i$ in the first position. As in the proof of Lemma~\ref{lem:chain}(c), pinching shows that the trace norm of this matrix is at least that of its leading $2\times2$ block. Since $\bigl\|\bigl[\begin{smallmatrix}a&c\\c&0\end{smallmatrix}\bigr]\bigr\|_1=\sqrt{a^2+4c^2}$ for real $a$ and $c$, this gives
\[\|s_ix_0E_{11}+Q_0\|_1\ge\sqrt{(48\mu_i)^2+4\cdot12^2}=24\sqrt{4\mu_i^2+1}.\]
Since $n_0\rho_0(x)=\sum_{i=1}^{12}\|s_ix_0E_{11}+Q_0\|_1$, the values of $\mu_i$ show that the term with $\ell=0$ contributes at least
\[24\bigl(\sqrt{145}+2\sqrt{37}+3\sqrt{17}+6\sqrt5\bigr)=1199.82\ldots\]
So
\[\frac{F_m(x)}{\|x\|_1}>\frac{136.5+1199.82}{667}=\frac{1336.32}{667}=2.00347\ldots>2.003.\qedhere\]
\end{proof}

\begin{remark}
The vector $u$ in the proof of Proposition~\ref{prop:cert} comes from a numerical maximiser. We evaluated $F_m$ using Theorem~\ref{thm:norm}(a) and maximised $F_m(x)/\|x\|_1$ as a function of $\sqrt x$ by a numerical quasi-Newton method. For $m=(40,\dots,40)$, the square roots of the resulting symmetric maximiser, scaled to have first entry $24$, are approximately
\[(24,\,6.41,\,4.89,\,3.63,\,2.56,\,1.63,\,0.80),\]
and $\sn{T_C(m)}\approx2.0043$. Rounding this vector gives $u$. For the vector $x$ of the proof, we have $F_m(x)/\|x\|_1\approx2.0037$.
\end{remark}

\begin{proof}[Proof of Theorem~\ref{thm:main}]
By Proposition~\ref{prop:cert}, the graph $C=(J_3-I_3)\otimes(J_4-I_4)$ is a productive core, so Theorem~\ref{thm:halfline} gives a real number $t$ with $[t,\infty)\subseteq\N$. The set $\N$ then has non-empty interior and is uncountable.
\end{proof}

\section{Numerical remarks and open questions}\label{sec:numerics}

In this final section, we make the number $t$ of Theorem~\ref{thm:main} explicit, compare the constants of the proof with numerical computations and record some open questions.

For the core of Proposition~\ref{prop:cert}, we have $\sn{C}=2$, $\sn{T_C((\overline{40})^{[6]})}>2.003$ and $\beta=5$, so $R_1=\max\{40,55\}=55$. Then $g>1.5\cdot10^{-3}$, and we take $R=55$, so that
\begin{multline*}\delta_{55}>\Bigl(1-\frac4{55}\Bigr)\frac{(1.5\cdot10^{-3})^2}{5^3\cdot55}>3\cdot10^{-10},\\ K\le\left\lceil\left(\frac{\log(4+3\cdot10^{-10})}{\log4}-1\right)^{-1}\right\rceil<1.85\cdot10^{10}.\end{multline*}
For $R=55$, the proof of Theorem~\ref{thm:halfline} gives the interval $[a,b]\subseteq\N$ with $a=\sn{T_C(\overline{55})}$ and $b=\sn{T_C(\overline{110})}$, where $2.003<a<b\le4$ and $b-a>3\cdot10^{-10}$. By Theorem~\ref{thm:halfline}, we have $[4^K,\infty)\subseteq\N$, and hence $[t,\infty)\subseteq\N$ for $t=4^{1.85\cdot10^{10}}$. Numerical experiments suggest that in fact $[t,\infty)\subseteq\N$ for considerably smaller, albeit still large, values of $t$. The size of $t$ is mainly due to the lower bound $\delta_R$ in Theorem~\ref{thm:large}(a) for $b-a$. We can only prove that $b-a>3\cdot10^{-10}$, while numerical computations give $b-a\approx3.9\cdot10^{-3}$. Even this length would only give $t\approx10^{109}$ by Lemma~\ref{lem:powers}, so the covering by powers of a single interval also limits the method.

To date, the core $C$ of Section~\ref{sec:proof} is the only core that we have proved to be productive. For $C=J_{15}-I_{15}$, numerical maximisation of $F_m$ gives $\sn{T_C(\overline{R}^{[12]})}>2$ at $R=1000$. Since $\sn{J_n-I_n}<2$ and $T_{J_{15}-I_{15}}(m)$ is an induced subgraph of $T_{J_n-I_n}(m)$ for $n\ge15$, this would make $J_n-I_n$ productive for every $n\ge15$. For $R=1000$, our numerics suggest that the norm is about $2.08$ for $n=500$, and for $n=100$ that it is already about $2.058$ at $R=55$, which would give a smaller $t$, but we have not certified these computations.

We conclude with several open questions about the structure of $\N$. Let $\eta_6=\frac25\sqrt{5+2\sqrt5}\approx1.231$. The set $\N\cap[0,\eta_6]$ is finite by~\cite[Theorem~1.1]{L14}, while $\N$ contains a half-line by Theorem~\ref{thm:main}, and Questions~\ref{q:astar}, \ref{q:two} and~\ref{q:blocky} concern how $\N$ passes between these two regimes. Let $a_*$ be the infimum of the numbers $a$ for which $[a,b]\subseteq\N$ for some $b>a$, and let $t_*$ be the infimum of the numbers $t$ with $[t,\infty)\subseteq\N$. Then $\eta_6\le a_*$ by the finiteness of $\N\cap[0,\eta_6]$. The proof of Theorem~\ref{thm:halfline} gives $a_*<4$, and the computation above gives $t_*\le4^{1.85\cdot10^{10}}$.

The number $2$ is a natural threshold for our method, and we do not know whether it plays any special role in the structure of $\N$. The decay of the weights in Proposition~\ref{prop:decay} needs norms above $2$, so the method only produces intervals above $2$. On the other hand, $2$ lies in $\N$, since it is the norm of the map on $\B(\ell^2)$ that takes an operator to its off-diagonal part, by~\cite[Theorem~2.1]{BCD} and Lemma~\ref{lem:finite}.

\begin{question}\label{q:astar}
\begin{enumerate}
\item[(a)] Is $a_*\le2$?
\item[(b)] Is $t_*\le2$, so that $[2,\infty)\subseteq\N$?
\end{enumerate}
\end{question}

There are points of $\N$ which, suggestively, lie only a little way above $2$. Indeed, the numbers $\sigma_k=\tfrac2{k+1}\cot\tfrac{\pi}{2k+2}$ are norms of paths by~\cite[Theorem~4.2]{L14}, so $\N$ contains $\sigma_5\sigma_9\sigma_{122}\approx2+1.3\cdot10^{-6}$ by Lemma~\ref{lem:tensor}. Consequently, any gap theorem such as $\N\cap(2,2+\varepsilon]=\emptyset$ could only hold for $\varepsilon<1.4\cdot10^{-6}$.

\begin{question}\label{q:two}
\begin{enumerate}
\item[(a)] Is $2$ a limit point of $\N$ from above?
\item[(b)] Is $[2,2+\varepsilon]\subseteq\N$ for some $\varepsilon>0$?
\end{enumerate}
\end{question}

Numerical experiments with smaller branching numbers suggest positive answers to Question~\ref{q:astar}(a) and Question~\ref{q:two}. For our core and the map $m\mapsto\sn{T_C(m)}$ on $\{10,11,\dots,20\}^{\bN}$, the inequalities~\eqref{eq:interval-hyp} appear to hold in every case that we tested, and $\sn{T_C(\overline{10})}$ appears to equal $2$. Together with continuity of this map, this would give $[2,b]\subseteq\N$ for some $b>2$. We have not certified these computations, the method of Section~\ref{sec:productive} does not apply to such small branching numbers. Moreover, we have no evidence concerning Question~\ref{q:astar}(b).

The next question combines~\cite[Question~4.8]{L14} with the remaining parts of~\cite[Question~4.9]{L14}. The limit points of $\N$ include $2$, $\sqrt2$ and $4/\pi$ by~\cite[Remark~4.7]{L14}.

\begin{question}\label{q:closed}
Is $\N$ closed? Is the infimum of the limit points of $\N$ equal to $4/\pi$, and are there limit points of $\N$ from above that are not limit points from below?
\end{question}

The graphs $T_C(m)$ have small block complexity in the sense of~\cite{BGHJN}. A $0$--$1$ matrix is \emph{blocky} if its support is a union $\bigcup_iX_i\times Y_i$ for families $\{X_i\}$ and $\{Y_i\}$ that are each pairwise disjoint, and by~\cite{Liv}, the non-zero blocky matrices are the $0$--$1$ matrices of norm $1$. For $L\ge1$, let $\N_L$ be the set of norms $\sn B$ of the $0$--$1$ Schur multipliers $B$ indexed by $\bN\times\bN$ with $B=\sum_{i=1}^L\sigma_iB_i$ for some blocky matrices $B_i$ and $\sigma_i\in\{-1,1\}$. Then $\N_1=\{0,1\}$ and $\N_L\subseteq[0,L]$, so no $\N_L$ contains a half-line. The set $\N$ is the union of the sets $\N_L$ by~\cite[Theorem~1.1]{BGHJN}, and $\N_K\N_L\subseteq\N_{KL}$ by Lemma~\ref{lem:tensor}, since the tensor product of two blocky matrices is blocky. For our core, we have
\[C=J_3\otimes J_4-J_3\otimes I_4-I_3\otimes J_4+I_3\otimes I_4,\]
and each of these four terms is blocky. Every other edge of $T_C(m)$ joins a vertex $v$ with $d(v)\ge1$ to its parent $v^-$. The edges $vv^-$ with $v^-$ a row, together with the edges of $I_3\otimes I_4$, have a blocky matrix, since the vertices of $C$ have no parents. The edges $vv^-$ with $v^-$ a column also have a blocky matrix. Hence every $\sn{T_C(m)}$ lies in $\N_5$, and the proof of Theorem~\ref{thm:halfline} gives $[a,b]\subseteq\N_5$ with $a<b$. So the least $L$ for which $\N_L$ has non-empty interior satisfies $2\le L\le5$. If $L=2$, then $a_*<2$, since $\N_2\subseteq[0,2]$.

\begin{question}\label{q:blocky}
What is the least $L$ for which $\N_L$ has non-empty interior?
\end{question}

\section*{Acknowledgements}

Anthropic's Claude large language models were used to explore proof strategies, draft text and write code for numerical computations.

\end{document}